\documentclass[11pt]{amsart}
\usepackage[margin=3cm]{geometry}
\title[A limit law for the sequential model of preferential attachment graphs]{A logical limit law for the sequential model of preferential attachment graphs} 
\author{Alperen \"{O}zdemir} 
\address{Department of Mathematics, KTH, Royal Institute of Technology} 
\email{alpereno@kth.se} 

\usepackage{amssymb, amsmath, float, graphicx,  mathtools, xcolor, mathdots}
\usepackage[colorlinks=true, urlcolor=black,linkcolor=blue, citecolor=blue]{hyperref}
\usepackage{enumerate}

\newtheorem{theorem}{Theorem}[section]
\newtheorem{lemma}{Lemma}[section]
\newtheorem{corollary}{Corollary}[section]
\newtheorem{definition}{Definition}[section]
\newtheorem{example}{Example}[section]

\newtheorem{remark}{Remark}[section]
\newtheorem{proposition}{Proposition}[section]

\newcommand{\bbox}{\hfill $\Box$}
\newcommand{\pf}{\noindent {\it Proof:} }

\newcommand{\floor}[1]{\left \lfloor #1 \right \rfloor}
\newcommand{\ceil}[1]{\left \lceil #1 \right \rceil}

\newcommand{\cal}[1]{\mathcal{#1}}
\newcommand{\tn}[1]{\textnormal{#1}}

\begin{document}

\begin{abstract}
For a sequence of random graphs, the limit law we refer to is the existence of a limiting probability of any graph property that can be expressed in terms of predicate logic. A zero-one limit law is shown by Shelah and Spencer for Erd\"{o}s-Renyi graphs given that the connection rate has an irrational exponent. We show a limit law for preferential attachment graphs which admit a P\'{o}lya urn representation. The two extreme cases of the parametric model, the uniform attachment graph and the sequential Barab\'{a}si-Albert model, are covered separately as they exhibit qualitative differences regarding the distribution of cycles of bounded length in the graph.
\end{abstract}

%, while one can find counterxamples with oscillating probabilities for rational exponents.
%see Section \ref{model} for the definitions,

\maketitle

\section{Introduction}

We can transcribe any graph into the language of predicate logic by taking our domain the set of vertices and associating edges with binary relations. Each such graph will define a model, for which we can ask if it satisfies a given logical sentence. We are interested in sequences of random graphs with an increasing number of vertices. We ask if the probability that a sentence is satisfied has a limit as the size of the random graph goes to infinity. 

A commonly raised question here is, if the limit exists, whether that limit is always zero or one for a fixed sentence. See the following surveys on zero-one logical law \cite{C89, W93, A18}. The zero-one law appears as a threshold phenomenon for parametric models, such as Erd\"{o}s-Renyi graphs $G(n,p_n)$. For example, the probability that there exists a cycle in $G(n,p_n)$ goes to zero if $n p_n\rightarrow 0,$ and goes to one if $n p_n \rightarrow \infty,$ which means that the threshold value is $n^{-1}.$

The first-order theory of random graphs is covered in detail in Spencer's book \cite{Sp01}. Not all graph properties can be expressed in terms of predicate logic, such as connectivity, Hamiltonicity, 2-colorability etc. require second-order terms. An interesting result is that, for those that can be expressed in the first-order logic, Erd\"{o}s-Renyi graphs satisfy a zero-one law for $p=n^{-\alpha}$ only if $\alpha$ is irrational \cite{SS88}. Furthermore, the authors constructed a sentence for which the limiting probability does not exist for a fixed rational exponent, and one can find such a sentence for any rational, see Section 8 of \cite{Sp01}. 

 In this paper, we study a model of sequentially constructed random graphs. The preferential attachment model was first suggested in \cite{BA99} as a model for real-world networks, with an emphasis on the power law distribution, which refers to that the ratio of vertices with degree $k$ decays proportional to $k^{-\alpha}$ for some $\alpha>0.$ The power law distribution has empirical support, going back to \cite{AJB99}, and challenges \cite{CSRNM09} as well. In the case of the uniform attachment model, the degree distribution shows a geometric law, which can be found in Section 4 of \cite{BRST01}. See also Durrett's book \cite{Du07}, Section 4, for several different approaches in the study of these graphs. One particular approach in defining that we will use is the sequential model, in which the edges are also attached sequentially along with the vertices. It was first defined in \cite{BBCS05} in order to make use of P\'{o}lya urns for finding bounds on the degrees of vertices.

Suppose we fix the connection probability of vertices in a sequence of Erd\"{o}s-Renyi graphs, and let the number of vertices goes to infinity. Then the limiting countably infinite graph is the Rado graph, see the survey \cite{C97}, which is unique up to isomorphism \cite{ER63}. For preferential attachment graphs, it is shown in \cite{KK05} that such unique limiting graphs exist for the cases where the fixed number of edges attached at every stage is either one or two.  If it is larger than two, one can find two non-isomorphic attachment graphs with infinite vertices. We are interested in a weaker notion than isomorphism, which is called elementarily equivalence and will be defined in Section \ref{elemequiv}. Now we state our result. 

\begin{theorem}\label{mainthm}
Let $G_n$ be a preferential attachment graph sampled according to the rule defined in \ref{defnt}. For any first-order sentence $\varphi$ over simple graphs,
 \[\lim_{n \rightarrow \infty}\mathbf{P}( G_n \vDash \varphi) \textnormal{ exists.} \] 
\end{theorem}
The result is likely to be extendible to the classical model as both models have  asymptotically the same weak local limit and the same degree distributions. However, it is hard to expect any strange result as in the case of Erd\"{o}s-Renyi graphs. The reason is that the local dynamics is much more stable in the attachment graphs, which play an important role in determining the limit law. In general, it could be interesting if there can be found any relation between the limit laws and the graph properties, such as the size of the diameter or denseness. For example, the diameter of a typical Erd\"{o}s-Renyi graph is of constant order \cite{CL01}, while it is almost of logarithmic order in the Barab\'{a}si-Albert model \cite{BR04}.

Earlier results on the limit law are obtained under a restriction on the length of the sentences \cite{M22a,MZ22b} or a restriction on the degrees of vertices \cite{MZ22}. In the latter, the authors use inhomogeneous processes on the neighborhoods of bounded depth to prove the convergence for bounded degree graphs. We use similar processes and their idea to speed them up. The main difficulty therein was to extend it to the neighborhoods of arbitrarily growing size.

In Section \ref{prelim}, following the logical preliminaries, we discuss the limitations of a straightforward approach, then define the prototypical processes that we will modify and use later in the paper. In Section \ref{prefattach}, we present the definitions and weak local limit results for the preferential attachment graph. In order to show convergence results using inhomogeneous processes, we first need to show that their transition probabilities converge individually. This part extensively relies on the weak local limit law in \cite{BBCS14}.

The following section is about the evolution of the neighborhoods of vertices and cycles. The distribution of the neighborhoods of cycles of bounded length plays a determinantal role in the limit law. The cycle generation probabilities are formulated using the P\'{o}lya representation of the model and the bounds on the those probabilities are obtained by using the ideas in \cite{BHJL23}, where they study the cycle density in the uniform attachment graphs. 

In the final and the most extensive section, we distinguish the absorbing logical classes and show that the fixed neighborhoods are attracted to there eventually, which resolves the problem with unbounded degrees of vertices.

% We later describe the method to show the logical limit for those graphs.
\section{Preliminaries} \label{prelim}
 
We provide the logical background in this section. We give some definitions, examples and introduce the tools to state and prove statements related to the convergence theorems mentioned in the introduction.

\subsection{Model theory}\label{model}
 We refer to \cite{CK90}, \cite{H93} and \cite{M06} for various different presentations of the model theory. First, we define a \textit{structure} $M$, through its four components:
\begin{enumerate}
\item A non-empty set $A,$ which is called the \textit{domain} of $M,$ 
\item A set of functions $F^A$ and positive integers $n_f$ such that $f^A:A^{n_f}\rightarrow A$ for $f \in F^A,$
\item A set of relations $R^A$ and positive integers $n_R$ such that $R^A \subseteq A^{n_R}$ for each $R \in R^A,$
\item A set of constant elements $C^A \subseteq A.$
\end{enumerate}
Any of the sets $F^A, R^A$ and $C^A$ can be empty.

We define a \textit{language}  $L$ to be a collection of non-logical symbols which include the symbols representing functions, relations and constants. We assume that the symbols can be read off from a given structure. On the other hand, the logical symbols are the negation ($\neg.$); the equality sign ($=$); the universal ($\forall$) and existential ($\exists$) quantifiers; and the Boolean connectives ($\vee, \wedge, \neg, \Rightarrow, \Leftrightarrow$).

%Taken from \cite{Sp01}.
%A \textit{signature} of A is a particular choice for all the components above.

We have \textit{variables}, which are symbols such as $x_1,x_2,\ldots.$ They substitute the elements of the structures. A variable is \textit{free} if it is not bound by any quantifier, that is to say  $\forall$ or $\exists$ does not appear in the formula. We define the set of $L$\textit{-terms} as the smallest set that contains
\begin{enumerate}[i)]
\item Every constant of $L,$
\item every variable $x_i$ for $i=1,2,\ldots,$
\item the expression $f(t_1,\ldots,t_{n_f})$ for every function  $f$ of $L$ and every set of terms $t_1,\ldots,t_{n_f}.$
\end{enumerate}
Then an \textit{atomic formula} is either 
\begin{enumerate}[i)]
\item $t_1=t_2$ if $t_1$ and $t_2$ are terms, or
\item the expression $R(t_1,\ldots,t_{n_R})$ if $R$ is a relation of $L$ and $t_1,\ldots,t_{n_R}$ are terms.
\end{enumerate}

Formulas are derived from atomic formulas by applying logical symbols listed above. A \textit{sentence} is a formula with no free variables, and a \textit{theory} is a set of sentences. For a structure $M,$ we say $M$ is a \textit{model} of a sentence $\varphi,$ if $\varphi$  is true in $M,$ denoted by $M  \vDash \varphi .$ If all sentences of a theory is satisfied by $M,$ then we write $M \vDash T.$  

 %On the existence of a model,

%A sentence is $atomic$ if it is made of a variable only. 
\begin{example} \normalfont 
Consider a simple graph $G=(V,E)$ with the vertex set $V$ and the edge set $E$. We can define the graph by a single binary relation. Let $V$ be the domain and the ordered pair of elements of $V$ lie in $R$ if there is an edge joining them. So for an undirected graph, we include both pairs in $R.$ Letting $L=\{\sim\},$ we have $u \sim v$ if and only if there is an edge connecting $u$ and $v.$ For example,
a graph is loopless if
\[\forall v \forall u \ ([u=v]  \Rightarrow   \neg [u \sim v]),\]
% $v \sim u$ if
and a graph is complete if
\[\forall v \forall u \ (\neg [u=v]  \Rightarrow  [u \sim v]).\]
\end{example}

The models we will study involve \textit{multigraphs}, the graphs with possibly more than one edge between two vertices. See \cite{K73} for the logic of multigraphs.

The last definition will be on the length of the sentences. For any formula $\varphi$ belonging to the first-order logic, the \textit{quantifier rank} of $\varphi,$ denoted by $\tn{qr}(\varphi),$ is inductively defined as:
\begin{enumerate}
\item If $\varphi$ is atomic, then $\tn{qr}(\varphi)=0,$
\item $\tn{qr}(\varphi)=\tn{qr}(\neg \varphi),$
\item $\tn{qr}(\bigwedge \Phi)=\tn{qr}(\bigvee \Phi)=\max \{ \tn{qr}(\psi):\psi \in \Phi \},$
\item $\tn{qr}(\forall x  \,\psi)=\tn{qr}(\exists x \, \psi)=\tn{qr}(\psi)+1.$
\end{enumerate}

\subsubsection{Elementary equivalence} \label{elemequiv}

We will classify the models in terms of the sentences that they satisfy. The following equivalence relations assume that the structures that are compared are defined on the same language.

\begin{definition}\label{eleqclass} \normalfont
The two models $\cal{A}$ and $\cal{B}$ are \textit{elementarily equivalent}, denoted by $\cal{A} \equiv \cal{B},$ if their truth values agree on all first-order sentences. They are called \textit{k-elementarily equivalent}, denoted by $\cal{A} \equiv_k \cal{B}$, if they agree on all sentences with quantifier rank less than or equal to $k$.
\end{definition}
Corollary 3.3.3 in \cite{H93} reads:

\begin{theorem} \label{finiteness}
For any two models $\cal{A}$ and $\cal{B},$ $\cal{A} \equiv \cal{B}$ if and only if $\cal{A} \equiv_k \cal{B}$ for all $k \in \mathbb{N}.$
\end{theorem}

Note that if two structures are isomorphic, that is to say if there is a bijection between two structures that preserves the relations and is compatible with functions, then they satisfy the same sentences. Therefore elementary equivalence is a weaker notion of similarity than isomorphism. See also Example 3.3.2 in \cite{Sp01} and Chapter 5 of \cite{H93} for more insight on this topic.

\begin{definition} \normalfont Let two structures $\cal{A}$ and $\cal{B}$ with a common language have domains $A$ and $B$ respectively. Provided that $S_A \subseteq A$ and $S_B \subseteq B,$ a function $g:S_A \rightarrow S_B$ is called a \textit{partial isomorphism} if it is a bijection that preserves all relations and functions of $L.$ 
\end{definition}

One of our concerns is the cardinality of the logical equivalence classes.

\begin{theorem}\label{cardlogic}
If a language $L$ contains only relations, there are finite number of equivalence classes. 
\end{theorem}
\noindent Observe that all possible relations over $k$ elements are finite. The theorem follows from a reverse induction argument carried out for binary words, Lemma 2.3 in \cite{L93}. In particular, the conclusion holds for the graphs, see Theorem 2.2.1 in \cite{Sp01}.  Although the number of equivalence classes does not depend on the size of the domain of the structures, it can be huge. For example, a lower bound on the number is given in  \cite{Sp01} by a tower function which is defined by the recursive formula $T(k)=2^{T(k-1)}.$

%In the case of groups, we cannot have such bound. Example, maybe?

\subsubsection{Ehrenfeucht-Fra\"{i}ss\'{e} Games}

To verify elementary equivalence between two structures, we use a perfect information, sequential and two-player game where the existence of a winning strategy for a player implies the elementary equivalence. We describe it below.

\begin{definition} \label{game}
The Ehrenfeucht-Fra\"{i}ss\'{e} game on two sets $A$ and $B$ with $k$ rounds, denoted by $\tn{EF}_k[A,B],$ is played between two players
\begin{align*}
\tn{Player  \textbf{I} } & \, \tn{ a.k.a \, Spoiler}, \\\tn{Player \textbf{II} }  & \, \tn{ a.k.a  \, Duplicator.}
\end{align*}
At each round, Player \tn{\textbf{I}}, first chooses one of the two sets $A$ and $B,$ then chooses an element from that set. Player \tn{\textbf{II}} responds by choosing an element from the other set. Let us denote the element chosen from the set $A,$ by any of the two players, at the $l$th stage by $a_l.$ We similarly define $b_l.$  The game is a \textit{win} for Player \tn{\textbf{II}} if there exists a partial isomorphism:
 \[g:\{a_1,\ldots,a_k\}\rightarrow \{b_1,\ldots,b_k\} \tn{ such that }g(a_i)=b_i \tn{ for all } i=1,\ldots,k. \]
\end{definition}
See Theorem 2.4.6 in \cite{M06} for the following: 
\begin{theorem}  \label{gameequiv}
The game $\tn{EF}_k[A,B]$ is a win for Player \tn{\textbf{II}} if and only if $A \equiv_k B.$
\end{theorem}
A worthwhile remark in \cite{Sp01} is that the advantage of the first player is to be able to alternate between the two structures.

\subsection{Inhomogeneous random processes}\label{inhomog}

%In all of the examples below, we will consider, so called growth processes. That is to say the number of elements in the structure increases by time. 

%We will review the random processes that we use in the paper. 
In this paper, we use various random processes to model the evolution of the neighborhoods in the graph for a fixed isomorphism class of a neighborhood. The reason for considering the neighborhoods will be clarified in Section \ref{evolution}. Those neighborhoods will consist of trees to a large extent but also of cycles which play a determinantal role for the limit. 

To highlight the difficulty, let us start with a naive approach and consider a Markov chain on the logical equivalence classes of the neighborhoods. Since there are finitely many equivalence classes, we will have a finite Markov chain. A simple observation will reveal that the probability that a cycle will be created is of decaying order, which will be shown to be $n^{-1}.$ A second fact, which is not hard to argue for, is that there are diverging order of trees of every feasible equivalence class. So we will have an inhomogenous Markov chain with a transition matrix converging to the identity matrix. This will not be good enough:

\begin{example} An inhomogeneous Markov chain with two states:\footnote{The author could not retrieve the reference for this example.} 
\[P_n=\begin{cases}
\begin{aligned}
\begin{pmatrix}
1-\frac{1}{n} & \frac{1}{n} \\
0 & 1
\end{pmatrix}  & \tn{ if } 2^{2k-1}<n<2^{2k}, \\[1ex]

\begin{pmatrix}
1 & 0 \\
\frac{1}{n} & 1-\frac{1}{n}
\end{pmatrix}  & \tn{ if } 2^{2k}<n<2^{2k+2}. 
\end{aligned}
\end{cases}\]
Observe that $P_n \rightarrow I$ but
 $\vec{x}\left(\prod_{i=1}^nA_i\right)$
  oscillates along the line between $\begin{pmatrix}
1 \\ 0
\end{pmatrix}$ and
$\begin{pmatrix}
0 \\1
\end{pmatrix}$ for any initial state $\vec{x}.$
\end{example}

 %will be considered separately. The number of two cycles has finite expectation, so we can ignore. Let us start with rare events.  Doeblin discussed this for non-homogeneous chains. 

%So we will consider the $k-1$ neighborhoods of cycles. If there are enough of them , $k$ many of them, the Duplciator can duplicate. Otherwise, the Spoiler can point out the elementarily difference between two graphs. 

However, if the Markov chain converges to an irreducible ergodic chain, the convergence to stationary distribution is guaranteed.
 
 \begin{definition} \normalfont \label{matrixnorm}
Let $A=(a_{ij})$ be a possibly countable square matrix. Define the norm  \[\|A \|_{\infty} = \sup_{i} \sum_{j} a_{ij}.\]
\end{definition}

\begin{theorem} \label{conv}  \tn{(Theorem 242 of \cite{V22})} 
Let $\{P_n\}_{n \in \mathbb{N}}$ be the sequence of transition matrices of a countable state space Markov chain. If we have $\|P_n-P\|_{\infty} \rightarrow 0$ for some ergodic stochastic matrix $P,$ then the inhomogeneous Markov chain is also ergodic.
\end{theorem}

 The asymptotic behavior of inhomogeneous chains attracted the attention of D\"{o}eblin \cite{D37}. For time-inhomogeneous Markov chains over a finite state space, D\"{o}eblin provided conditions on time-dependent transition probabilities for the decomposition of the state space into asymptotically closed communicating classes. In \cite{B45}, with no further assumption than the finiteness of the state space non-homogeneous chains, Blackwell proved that the space-time representation of the chain can be decomposed into  finite number of  sequences such that the chain will eventually remain in exactly one of them with probability one. Even under certain regularity conditions, it is not much known about the precise distribution and the rate of convergence of inhomogeneous chains, see \cite{SZ07} for a spectral analysis of inhomogeneous chains. But since we are interested only in the existence of limit, those chains turn out to be useful in this context. 

 A particular setting we will use inhomogeneous processes is the set of isomorphism classes of neighborhoods of fixed depth, which admits a lattice structure in agreement with the vertex attachment. They will be defined in the next chapter. Here, we define the processes in a general sense and refer to them later.

\subsubsection{A martingale and a law of the iterated logarithm}\label{martaatl}

Consider the sequence of random variables $\{M_n\}_{n \in \mathbb{N}}$ defined inductively as
\begin{equation*} 
\begin{aligned}
   &   
  M_{n+1} =M_n + \begin{cases}   
       1,&  \text{with prob. } p(n),\\[5pt]
      -m,&  \text{with prob. } \frac{K M_n}{n+1}, \\[5pt]
       0,&   \tn{ otherwise}.
    \end{cases}
\end{aligned}
\end{equation*}
We assume that $\lim_{n \rightarrow \infty} p(n)$ exists, and $K$ and $m$ are positive integers. Letting $s=K m,$ we have
\[\mathbf{E}[M_{n+1} | \mathcal{F}_n]= \left( 1-\frac{s}{n+1}\right)M_n +p(n)\]
where $\cal{F}_n$ is the $\sigma$-algebra generated by $M_1,\ldots,M_n.$ Using the notation $(n)_k:=n (n-1)\cdots (n-k+1),$ we obtain the martingale
\[Z_{n+1}=(n+1)_s M_{n+1} - \mu_{n+1}\]
where its mean is recursively defined as
\[\mu_{n+1}=\mu_n+ (n+1)_s p(n).\]

The martingale differences are 
\begin{equation*} 
\begin{aligned}
   &   
  X_{n+1}:=Z_{n+1}-Z_n \, |\mathcal{F}_n = s (n)_{s-1} M_n - (n+1)_s p(n)+ \begin{cases}   
       (n+1)_s,&  \text{with prob. } p(n),\\[5pt]
      -m(n+1)_s,&  \text{with prob. } \frac{KM_n}{n+1},\\[5pt]
       0,&   \tn{ otherwise}.
    \end{cases}
\end{aligned}
\end{equation*}
We  note that since $\mathbf{E}[X_n]=0,$ we have 
\[\mathbf{E}M_n=\frac{(n+1)}{s}\mathbf{E}p(n)\]
Then, we can show that
\[\textnormal{Var}(X_{i+1})=(i+1)_s^2 p(i) + s(i+1)_s^2\frac{M_i}{i+1} \sim (i+1)^{2s-1}\mathbf{E}M_i. \]
where ``$\sim$" refers to the asymptotic equality. Suppose $p(n)\sim n^{-p}.$ Then, the variance of $Z_n$ is of order 
\[\sum_{i=1}^n (i+1)^{2s-1}\mathbf{E}M_i \sim n^{2s+1-p}.\]

%There are refined versions available, but we can take the following 
%\begin{theorem}
%Let $Z_n$ be a martingale with martingale differences $X_n.$ 
%\end{theorem}
Among many other possible probabilistic inferences we can derive from this martingale,  our particular interest is in the law of iterated logarithm. The following result is based on a slightly weaker result of Theorem 5.4.1 in \cite{S74}.
\begin{lemma} \tn{(Lemma 1 of }\cite{F92}\tn{)}\label{lil}
Let $\{Z_n, \cal{F}_n, n\geq 1\}$ be a martingale with martingale differences $X_i=Z_{i+1}-Z_i.$ Let $s_n^2$ be the variance of $Z_n$ such that $s_n^2 \rightarrow \infty.$ If there exists $\cal{F}_n$ measurable $K_n > 0 $ for all $n\geq 1$ such that  $ \sup_n K_n \leq K$ for some constant $K>0,$ and
\[|X_i| \leq K_i \frac{s_n}{\sqrt{2 \log_2 s_n^2 }} \quad \tn{a.s.},\]
then
\[\limsup_n \frac{Z_n}{2 s_n \sqrt{2 \log_2 s_n^2}} \leq C(K)\]
for some constant $C(K).$ 
\end{lemma}
We will use martingales in Section \ref{trees} and \ref{sabg}.

% We can isolate the $\frac{1}{n}$ term as it happens infinitely often.

\subsubsection{A slow Markov chain} \label{inhomMC}
We look at a special case of the martingale above, where we can identify the stationary distribution. Define
\begin{equation*} 
\begin{aligned}
   &   
  M_{n+1} =M_n + \begin{cases}   
       1,&  \text{with prob. } \frac{\rho(n)}{n+1},\\[5pt]
      -1,&  \text{with prob. } \frac{\tau(n) M_n}{n+1}, \\[5pt]
       0,&   \tn{ otherwise}.
    \end{cases}
\end{aligned}
\end{equation*}

Observe that the $M_{n+1}-M_n$ is non-trivial infinitely often as $\sum \frac{1}{n+1}=\infty.$ To obtain a homogeneos limit for the process, we want to eliminate the divisor $\frac{1}{n+1}$ from the transition probabilities. S we condition on $\{ M_n \neq 0 \},$ and define a new Markov chain with transition probabilities 
 \begin{equation*} 
W_n(i,i-1)=\frac{\tau(n)i}{\tau(n)i+\rho(n)}, \, W_n(i,i+1)=\frac{\rho(n)}{\tau(n)i+\rho(n)} \ \tn{ for }i\geq 1 
\end{equation*}
and $W_n(0,1)=1,$ otherwise zero. Now suppose $\rho(n)\rightarrow \rho$ and $\tau(N)\rightarrow \tau.$ By Theorem 242 in \cite{V22}, we have 
\[ \lim_{n \rightarrow \infty} \prod_n W_n = \Pi,\]
whose column entries that can be determined by solving the recurrence relation:
\begin{align*}
\pi_0&=\frac{1}{1+\rho/\tau} \ \ \tn{ and } \ \ \pi_n = \frac{n+\rho/\tau}{n}\pi_{n-1} - \frac{\rho/\tau (n+\rho/\tau)}{n(n-2+\rho/\tau)}\pi_{n-2} \ \tn{ for }n\geq 1.
\end{align*}

This recurrence relation has the solution
\[\pi_n=\frac{(n+\rho/\tau) (\rho/\tau)^{n-1}}{2e^{\rho/\tau}n!} \ \tn{ for }n\geq 0.\]
This gives the stationary distribution for the inhomogeneous chain, which is just a tilted Poisson distribution. We will couple this Markov chain and others defined on a more particular space with more elaborate transition rules in Section \ref{lattice} and onwards.

\section{Preferential attachment graphs} \label{prefattach}

In this section, we define the random graph model that we will study and list some facts related to it, in particular related to its local limit. We are considering a sequence of random graphs where at each stage a vertex is generated and is attached to already present vertices according to a rule which favors vertices with large degrees. Note that ``the sequential model" does not refer to this sequential construction. Formally, we define

\begin{definition} \label{defnt} \normalfont (Preferential attachment model) 
Let $\alpha \in [0,1]$ and $m$ be a positive integer. We construct a sequence of undirected multigraphs 
\[G_{1} \subset G_{2} \subset \cdots \subset G_{n} \subset \cdots  \]
%\[G_0 \subset \tn{PA}_{1} \subset \tn{PA}_{2} \subset \cdots \subset \tn{PA}_{n} \subset \cdots \]
where ``$\subset$" implies both the vertex and the edge set inclusion. Let $G_1$ be a vertex with no edge and the vertex set of $G_{n}$ be $[n]:=\{1,2,\ldots, n\},$ labelled by natural numbers. $G_{n+1}$ is constructed from $G_n$ by connecting $m$ vertices $v_1,\ldots,v_m,$ not necessarily distinct, to the vertex $n+1$ according to the rules of the models provided below.

\noindent \textit{The classical model:} We choose $m$ vertices $v_1,\ldots,v_m$ independently according to the following rule. With probability $\alpha$ we uniformly choose a  vertex from $[n]$ and connect it to $(n+1).$ With probability $1-\alpha,$ we choose in proportion to the degrees of vertices, that is, 
\[\mathbf{P}(v_i=k)=\frac{\tn{deg}_n(k)}{2m(n-1)}\]
where $\tn{deg}_n(k)$ denotes the degree of $k$ and $2m(n-1)$ is the sum of degrees of all vertices.

\noindent \textit{The sequential model:} In the classical model, we attach all edges at once. In the sequential model, we attach each one of the $m$ edges one-by-one and update the attachment probabilities after each assignment of edges. 
We take the new uniform attachment coefficient 
\begin{equation} \label{alpha}
\alpha_n(i)=\alpha \frac{2m(n-1)}{2m(n-2)+2m\alpha +(1-\alpha)(i-1)}\in \left[\alpha,\alpha+ \frac{1}{n-2}\right).
\end{equation}
Then letting
\[\tn{deg}_n(k,i):=\tn{deg}_n(k)+|1 \leq j \leq i-1 \, | \, v_j=k)|,\] we have instead
\[\mathbf{P}(v_i=k)=\frac{\tn{deg}_n(k,i)}{2m(n-2)+i-1}\]
with probability $1-\alpha_n(i).$ 
\end{definition}

The boundary case where $\alpha=1$ is called the \textit{uniform attachment model}, which agrees for both models. We will call the other boundary case, $\alpha=0,$ the \textit{sequential Barab\'{a}si-Albert model}.

%The authors in \cite{BBCS14} show the weak local limit of the preferential attachment model using this model, then a coupling argument, see Section 4 of the same paper, extends the the limit to the classical model and a variation of it where multiple edges are not allowed.

The main advantage of this model is that the degree of a vertex relative to the sum of degrees of vertices with lower indices can be counted with balls in P\'{o}lya urns. Then applying de Finetti's theorem one gets conditionally independent attachment probabilities, which are given by $\beta$-random variables \eqref{psi}. See also \cite{PRR17} and \cite{GHHR22} for other applications of this model.

\subsection{Weak local limit} \label{polya}

In this section, we will study fixed depth neighborhoods of a randomly chosen vertex from the graph by using the \textit{weak local limit} of the preferential attachment graphs, which was found in \cite{BBCS14}, and named P\'{o}lya point-graph. This notion is first defined by Benjamini and Schramm in \cite{BS01}, which is called the \textit{distributional limit} there. See \cite{AS04} and \cite{AL07} also.

Let $d(v,w)$ the length of the shortest path between $v$ and $w.$ We define the $r$-neighborhood of vertex a $v$ as 
\begin{equation} \label{distance}
B_r(G,v)=\{w \in G \, : \, d(v,w)\leq r\}.
\end{equation}

 If there is no ambiguity we will drop $G$ in the notation and simply write $B_r(v).$ The notation $G \cong H$ means that $G$ is isomorphic to $H.$

\begin{definition} \label{PAmodel}  \normalfont (Weak local limit)
Let $G_n$ be a sequence of finite graphs and $v_{n}$ be a uniformly chosen vertex from $G_n.$ We call a rooted random graph $(G,x)$ the weak local limit of $(G_n,v_n)$ if  
   \[\lim_{n \rightarrow \infty}\textbf{P}(B_r(v_n) \cong (H,y)) = \textbf{P}(B_r(x) \cong (H,y))\]
   for all finite rooted graphs $(H,y)$ and for all $r \in \mathbb{N}.$
 \end{definition}

% $x$ uniformly chosen from  $[0,1]$ such that 

%Next we define it at least to adapt its use to 

%Reproducing the argument in \cite{BBCS14}, we l
%variation and an interpolation with uniform attachment model.

\subsubsection{Sampling from the sequential model} \label{seqmodel} We let
\begin{equation} \label{uchi}
u:=\frac{\alpha}{1-\alpha} \tn{ and } \chi:=\frac{1+2u}{2+2u}
\end{equation}
 where $\alpha$ is defined in Definition \ref{PAmodel}. The urn process interpretation gives that the new vertex is attached to vertex $i,$ conditioned on it is attached to one of the vertices in $\{1,2,\ldots, i\},$ is given by the $\beta$-random variable
\begin{equation} \label{psi}
\psi_i \sim \beta(m+2mu, (2i-3)m+2mu(i-1))
\end{equation}
where ``$=_d$" denotes the equality in distribution. So, the probability that the new vertex is attached to vertex $i$ is 
\[\varphi_i := \psi_i \prod_{j=i+1}^n (1-\psi_j).\]
We also define
\[S_l=\sum_{i=1}^l \varphi_i.\]
For all $l=1,\ldots,n,$ take $\{U_{l,i}\}_{i=1}^m$ i.i.d. random variables that are uniform over $[0,S_{l-1}].$ We connect $j<l$ to $l$ if $U_{l,i} \in [S_{j-1},S_j)$ for some $i=1,\ldots, m.$ The resulting graph has the same distribution with $G_n.$ In particular, this gives us a procedure of sampling trees from $G_n.$

% we can sample a tree from the sequential model by choosing $m$ random points from $[0,S_{l-1}]$ for all $i$ and drawing an edge between $j$ and $l$ for all points in $[S_{j-1},S_j)$ among the $m$ points. 

\subsubsection{P\'{o}lya-point graph} \label{polyatree} On the other hand, the P\'{o}lya-point graph, a random tree sampled from $G_n,$ is defined recursively in \cite{BBCS14} as follows. The root of the neighborhood is chosen according to
\begin{equation} \label{root}
X_0=Y_0^{\chi} \tn{ where } Y_0 \sim \tn{Unif}[0,1].
\end{equation}
Then we consider the $m$ neighbors of $X_0$ which are generated before $X_0.$ Let us call them $X_{0,1},\ldots, X_{0,m},$ which are chosen uniformly from $[0,X_0].$  The remaining neighbors of the root is distributed according to a Poisson point process on the interval $[X_0,1],$ find it below \eqref{pol1}. In particular, there are randomly many of those points. Continuing this process, we denote a vertex with distance $s$ to the root by $X_{\bar{a}}$ where $\bar{a}$ stands for $a_1,\ldots,a_{s-1}, a_s,$ which is called the \textit{height} of the sampled point in \cite{BBCS14}.
%where  $a=(0,a_1,a_2,\ldots,a_s)$ is the sequence of the indices of its ascendants and itself. 
Then the vertices attached to $X_{\bar{a}}$ and which are to the left of it, are distributed as
\begin{equation*}
X_{\bar{a},i} \sim \tn{Unif}[0,X_{a}] \, \tn{ for }i=1,2,\ldots,m
\end{equation*} 
and $X_{\bar{a},j}$ for $j=m+1,\ldots, N$ are generated by Poisson point process with intensity 
\begin{equation}\label{pol1}
 \rho(x)= \frac{1-\chi}{\chi} \, \Gamma(m+\mathbf{1}_{\{a_s\leq m\}},1) \,   X_{\bar{a}}^{-\frac{1-\chi}{\chi}} \, x^{\frac{1-2 \chi}{\chi}} dx \quad \tn{ on } \ [X_{\bar{a}},1]
 \end{equation}
where $\mathbf{1}_A$ is the indicator function of a given event $A.$ The $\Gamma$-random variable shows up as a scaling limit of $\beta$-distribution, such as $ n \beta(k,n)\rightarrow_d \Gamma(k,1).$  This process gives a random rooted graph $T$ of any desired depth.
 
 The weak limit is proved by coupling the neighborhoods of fixed depth of the graph sampled from the P\'{o}lya urn model in (\ref{seqmodel}) and the P\'{o}lya-point tree $T$ in (\ref{polyatree}). In particular, they define an isomorphism $\mathbf{k}$ between a ball of radius $r$ in $T$ and a ball of the same radius in $G_n$ with roots $X_0$ in \eqref{root} and $\mathbf{k}(0)=\ceil{nX_0^{1/\chi}}.$ The latter indeed agrees with the the uniform distribution over $\{1,\ldots,n\}$. For example, one gets
\[\left|S_{\mathbf{k}(\bar{a})}-X_{\bar{a}} \right| \leq \varepsilon\]
with high probability for large enough $n$.

 We can make particular choices for the root in $T.$ Most often, we will take $X_0=1,$ which will be associated with the new vertex attached to the attachment graph.
 
 \begin{lemma}\label{coupling}
 There exist a coupling between the P\'{o}lya-point graph $(G,x_0)$ for a deterministic $x_0 \in [0,1]$ and the rooted sequential model $(G_n, v_n)$ where $v_n=\ceil{n x_0^{1/\chi}},$ which satisfies the following. For any $\varepsilon >0,$ there exists $N(\varepsilon)$ such that $n>N(\varepsilon)$ implies
\[\mathbf{P}(B_r(v_n) \cong B_r(x_0)) \geq 1-\varepsilon.\]
\end{lemma}

 \pf It suffices to follow the proof for the local convergence in \cite{BBCS14}, which is by induction on the depth of the tree constructed. For the basis step, they look at the neighborhood of the random root (Lemma 3.5),  and for the inductive step they extend the same idea to $r$-neighborhoods (Lemma 3.6).  As already noted in the proof of the latter, the former result applies to the neighborhood of any given vertex. For the uniform attachment graph, we  simply take $\chi=1.$
 
 \bbox
  
%In fact, one can Combining Lemma 3.5 applied to a determinisitic vertex and the coupling between sequential and independent models in Section 4 give the result.
%This is not generally true. The idea does not apply to disconnected graphs, for instance. The local limit can still be defined there. Refer to that survey...\\

\begin{remark}
\normalfont 
The coupling for the uniform attachment graph can be viewed as a degenerate case of the one defined above. We simply take $\varphi_i=\frac{1}{n}$ and $S_i=\frac{i}{n}.$ Suppose the randomly chosen root is $k_0$. Let us define $y_0$ by setting $k_0=\ceil{ny_0}$ and take $p_l=\frac{1}{l-1}.$ We let $X_{l,i}$ be a Bernouilli random variable with parameter $p_l.$ Then the number of vertices attached to $k_0$ which lie between $k_0$ and $\ceil{ny}$ for $y_0 < y \leq 1$ is given by
\[N_{y_0}(y)= \sum_{l=k_0}^{\ceil{ny}} \sum_{i=1}^{m} X_{l,i}=\sum_{l=k_0}^{\ceil{ny}} \tn{Binom}(m,p_l).\]
The same count in the associated P\'{o}lya-point graph is given by following inhomogeneous Poisson point process $\mathbf{N},$ defined as 
\begin{equation*} 
 \mathbf{P}(\mathbf{N}((y_0,y]=s))=\frac{\Lambda(y_0,y)^s}{s!}e^{-\Lambda(y_0,y)}
\end{equation*}
where
\begin{equation*}
\Lambda(a,b)=m \log \frac{b}{a}.
\end{equation*}
\end{remark}

\subsection{Degree distribution}

We will list some facts on the degree distributions of vertices, some of which can be found in Section 5 of \cite{BBCS14}. For a fixed vertex $k \in \{1,2,\ldots,a\}$ where $a\leq n,$ we define the truncated probabilities as
\[\varphi_k^{(a)}= \psi_k \prod_{i=k+1}^a (1-\psi_i)\]
where $\psi_i$ is defined as in \eqref{psi}. Observe that
\[\sum_{i=1}^a \varphi^{(a)}_i=1\]
for any $a \leq n.$ Next we define the intervals
\[I_{k}^{(a)}=\left[ \sum_{j=1}^{k-1} \varphi_j^{(a)}, \sum_{j=1}^{k} \varphi_j^{(a)} \right]=\left[S_{k-1}^{(a)},S_{k}^{(a)}\right]=\frac{1}{S_a} [S_{k-1},S_{k}].\]
Let us define a set of i.i.d. uniform random variables on $[0,1],$
\[\left\{ U^{(a)}_i \, \big| \, i=1,\ldots,m \, ; a=1,\ldots,n \right\},\]
and based on those we define
\[\chi^{(a)}_{k,i}=\begin{cases} 
      1 & \tn{ if } U^{(a)}_i \in I_{k}^{(a)} \\
      0 & \tn{ otherwise.}
   \end{cases}\]
  %Observe that the event $\left\{ U^{(b)}_i \in I_{a}^{(b)} \right\}$ is equivalent to  $\{S_{a-1}\leq U < S_a \}$ for $U$ uniform on $[0,1].$

Let $\mathcal{F}$ be the $\sigma$-algebra of $\{\psi_i\}_{i=1}^{\infty}.$  We have 
\begin{equation} \label{chidef}
\mathbf{E}\left[\chi^{(a)}_{k,i} \, | \, \cal{F}\right]= \varphi_k^{(a)}=\psi_k \prod_{i=k+1}^a (1-\psi_i).
\end{equation}

Observe that for any vertex $a$ and $i=1,\ldots,m,$ there is a unique vertex $k$ to which $a$ is attached, but $a$ can be attached to the same vertex for $i_1\neq i_2$. So, we have
\begin{equation} \label{chi}
\begin{aligned}
   &   
   \mathbf{E} \left[\chi_{k_1,i_1}^{(a_1)} \chi_{k_2,i_2}^{(a_2)} \, | \, \cal{F} \right] = \begin{cases}   
       0,&  \text{if } a_1=a_2, i_1=i_2 \text{ and } k_1\neq k_2  \\[5pt]
      \mathbf{E} \left[\chi_{k_1,i_1}^{(a_1)} \, | \, \cal{F} \right],&   \text{if } a_1=a_2, i_1=i_2 \text{ and } k_1 = k_2 , \\[5pt]
       \mathbf{E} \left[\chi_{k_1,i_1}^{(a_1)} \, | \, \cal{F} \right] \mathbf{E} \left[\chi_{k_2,i_2}^{(a_2)} \, | \, \cal{F} \right] ,&  \tn{ otherwise}.
    \end{cases}
\end{aligned}
\end{equation}

%If we write it more explicitly,
%\begin{align*}
%\mathbf{E} \left[\chi_{k_1,i_1}^{(a_1)} \chi_{k_2,i_2}^{(a_2)} \, | \, \cal{F} \right] & = \psi_{k_1} \psi_{k_2} \prod_{j_1=k_1+1}^{a_1} (1-\psi_{j_1}) \prod_{j_2=k_2+1}^{a_2} (1-\psi_{j_2})
%\end{align*}

 Let the random variable  $D_n(k,i)$ stand for the degree of $k$ after the attachment of the $i$th edge with the addition of the $n$th vertex. Let us set $D_n(k):=D_n(k,m).$ We have the decomposition
  \begin{equation} \label{degreedist}
D_n(k)= m+ \sum_{l=k}^{n-1}\sum_{i=1}^m\chi_{k,i}^{(l)}.
 \end{equation}
Thus,
  \begin{equation} \label{conddegreedist}
\mathbf{E}\left[D_n(k) | \cal{F}\right]= m+ m \sum_{l=k}^{n-1} \varphi_{k}^{(l)}.
 \end{equation}
The degree function scaled by $n^{1-\chi}$ is shown to have a limiting distribution in Lemma 5.1 of \cite{BBCS14}.

%, which is uniform in $n$. See Lemma 5.1 inSection 5 of \cite{BBCS14} for more statistical properties of it.

Let us make a few observations on the correlations of random degrees. First, observe that, 
\begin{equation} \label{sametarget}
\mathbf{E}[D_{n_1}(k) D_{n_2}(k)] \geq \mathbf{E}[D_{n_1}(k)] \mathbf{E}[D_{n_2}(k)]
\end{equation}
by the very attachment rule that favors degree. On the other hand,
\begin{equation} \label{difftarget}
\mathbf{E}[D_{n_1}(k_1) D_{n_2}(k_2)] \leq \mathbf{E}[D_{n_1}(k_1)] \mathbf{E}[D_{n_2}(k_2)]
\end{equation}
if $k_1\neq k_2$ because any edge connected to $k_1$ is not connected to $k_2.$ We will show that they are asymptotically uncorrelated.

\subsection{Evolution of the logical equivalence classes} \label{evolution}

To show Theorem \ref{mainthm}, we will use the Ehrenfeucht-Fra\"{i}ss\'{e} game defined in Section \ref{game} and combine it with Theorem \ref{finiteness}. The game would be relevant only if the vertices chosen in either of the two graphs can be connected by a path, which is to be shown in this section. Since it has also finite rounds, we only need to consider local structure of the graph, which makes the local limit of \cite{BBCS14} useful here. In fact, it will suffice to consider large enough neighborhoods in the graph. We will mainly focus on three different kind of neighborhoods; rooted trees, cycles, and multi-cycles. We will go over them one-by-one, to show their abundance, stationary behavior and rarity respectively. The case with the cycles will be the most challenging. Let us define the neighborhoods and establish the logical connection below. 

%Defining  In particular, our analysis will boil down to the study of the neighborhoods of vertices around cycles of a bounded size in a similar way as employed in \cite{MZ22}.   

%Let $d(u,v)$ the length of the shortest path between $u$ and $v.$ Then we define the $r$-neighborhood of vertex a $v$ as \[U_r(v)=\{u \in G \, : \, d(u,v)\leq r\}.\] 
Similar to \eqref{distance}, we define the distance of a vertex to a set of edges as the minimal distance of the vertex to any of the edges in the set. We define the $r$-neighborhood of a cycle $\tn{C}=\{v_1,\ldots,v_l\}$ as 
\begin{equation} \label{cycneigh}
B_r(\tn{C})= \{ u \in G \, :\, d(u,v) \leq r \tn{ for some } v \in \tn{C} \}.
\end{equation}
We refer to these neighborhoods as \textit{cycle-trees}.  

Next, we define a cycle profile of a graph. Let $\mathbf{C}_{2r}$ be the set of cycles of length at most $2r$ in $G_n.$ Then, let us define an equivalence relation on cycles as $\tn{C} \leftrightarrow \tn{C}'$ if there exists $\tn{C}_1,\ldots,\tn{C}_s \in \mathbf{C}_{2r}$ such that $\tn{C}_0=\tn{C},\tn{C}_s=\tn{C}'$ and $d(\tn{C}_{i-1},\tn{C}_i)\leq r$ for all $i=1,\ldots s.$ We let 
\[\cal{C}=\left\{ \bigcup\limits_{\tn{C} \in [C]} \tn{C} \, : \, [\tn{C}] \in \mathbf{C}_{2r} / \leftrightarrow \right \}.\]
If we need to specify the stage of the attachment graph, we will write $\cal{C}_r(G_n).$ In addition, we have the set of  \textit{isolated} cycles,
\[\cal{C}_1=\{\tn{C} \in \mathbf{C}_{2r} \,: \, d(\tn{C},\tn{C}') > r \tn{ for all }\tn{C}' \in \mathbf{C}_{2r} \tn{ such that } \tn{C}'\neq \tn{C} \},\]
which are equivalence classes of size one in the relation defined above. So $\cal{C} \setminus \cal{C}_1$ is the set of \textit{multicycles}.
Finally, we define the bounded neighborhoods of points that do not contain cycles as
\begin{equation}\label{nocycle}
\cal{C}_0=\left\{ v \in G_n\, : \, B_r(v) \cap \tn{C}=\emptyset \tn{ for all } \tn{C} \in \mathbf{C}_{2r} \right\}
\end{equation}

Now we look at the logical classes. Let $\cal{L}$ be the set of all logical classes, which is finite by Theorem \ref{finiteness}. For a simple graph $G$, we define the vector 
\[\Lambda(G)=\left( \Lambda_1(G),\ldots, \Lambda_{|\cal{L}|}(G)  \right)\]
 on $\mathbb{Z}_{\geq 0}^{|\cal{L}|}$ by letting
\[\Lambda_L(G) = \left| S \in \cal{C}(G) \, : \, S \tn{ has logical type }  L \tn{ as an induced graph of }G\right|\]  
 for all $L \in \cal{L}.$ 

We say that $G$ and $G'$ have the \textit{same cycle profile} if and only if $\Lambda(G)=\Lambda(G').$ Then, we say that they \textit{agree on acyclic neighborhoods}, if the same vector we define over trees instead of cycles agree on both graphs. However, since it will be shown that the representatives for every feasible class is arbitrarily large when it comes to tree neighborhoods, we refrain from additional notation.

%contingent on the two possible cases of the choice of , let us call it $v.$ 

\begin{lemma}\label{logic}
If $G$ and $G$ have the same cycle profile, and they both have at least $k$ distinct neighborhoods for every logical class associated with the rooted trees, then $G \equiv_k G'$.
\end{lemma} 

\pf 
Following the strategy in Section 2.4 and 2.6.3 of \cite{Sp01}, we consider neighborhoods shrinking in terms of depth at each round. We start with the depth $r=\frac{3^k+1}{2},$ which is obtained from the recursion $r_1=1$ and $r_{i+1}=3r_i+1$ as $r_k.$ Let us define the inverse sequence $d_i:=r_{k-i-1}$ for notational convenience. This choice of large enough coefficient for the recursion is shown to clear the ambiguity for Player \textbf{II} whether a new vertex chosen by Player \textbf{I} is to be considered in the neighborhoods of earlier ones or independent from them. 

Let us describe a winning strategy for Player \textbf{II} under the hypothesis. At the first round, whichever vertex of any of the two graphs Player \textbf{I} chooses, Player \textbf{II} can always find a vertex in the other graph whose neighborhood share the same logical class because of the agreement of neighborhoods of depth $r$, cyclic or acyclic.

Suppose $v_1,\ldots,v_{i-1} \in G$ and $w_1,\ldots,w_{i-1} \in G'$ are the vertices chosen by the $i$th round. Assume Player \textbf{I} chooses $v_i \in G$ without loss of generality.  If $v_i$ belongs to a $2d_i-$neighborhood of any of previously chosen vertices in $G$, say $v_j,$ then Player \textbf{II} has a choice in the neighborhood of the associated vertex $w_j,$ because their neighborhoods of depth $2d_j$ logically agree. However, if $v_i$ is outside of that range for any of the $(i-1)$ neighborhoods of earlier vertices, then it can be treated as independent from them. 

Suppose $d(v_i,v_j)\geq 2d_j+1$ for all $j=1,\ldots,i-1.$ If $B_{d_i}(v)$ contains no cycles, then it can be duplicated because $G$ and $G'$ agree on acyclic neighborhoods. Suppose its $d_i$-neighborhood contains a cycle or multi-cycles, that is $S \subseteq B_{d_i}(v)$ for some $S \in \cal{C}(G)$. We do not rule out the case that $v_j \in B_{r}(S)$ for some $j < i.$ By the assumption that the graphs have the same cycle profile, there exists $\tn{S}' \in \cal{C}(G')$ such that $U_{r}(S) \equiv_k B_{r}(S').$ This implies that Player \textbf{II} can find $w_j \in B_{d_i}(S')$ which will extend the partial isomorphism to a one between $v_1,\ldots,v_i$ and $w_1,\ldots,w_i.$ 
 
\bbox

In the following two chapters, we will show that the acyclic neighborhoods have enough many representatives and the cycle profile has a limiting distribution.

\section{Evolution of the neighborhoods}

\subsection{Trees} \label{trees}

 We are interested in the multiplicity of each tree in the graph up to isomorphism. We want to show that all trees have enough copies to verify the hypothesis of Lemma \ref{logic}. We say an isomorphism class of a tree is \textit{feasible} if it has positive probability to appear in the attachment graph. For instance any tree with a vertex of degree less than $m$ is not feasible.   
 
\begin{lemma}\label{treelimit}
Let $T$ be a feasible isomorphism class of trees with depth $r,$ and $T_n$ denote the number of copies of $T$ in $G_n$ as an induced graph. We have
\[\lim_{n\rightarrow \infty}\mathbf{P}(T_n\geq k)=1.\]
\end{lemma}
\pf Assume that the degree of the root of $T$ is $m,$ so that it can be generated by taking $n$ as its root. We will set up a martingale $M_n$ as in Section \ref{martaatl} to find a lower bound for $T_n$.  

Let $p_T(n)$ be the probability of generating such tree. By Lemma \ref{coupling}, this converges to the probability of generating $T$ in the P\'{o}lya-point graph, which implies  $p_T(n)\rightarrow p_T>0$  This immediately suggests $T_n$ is of order $n.$ A problem here is that the generation and the removal of trees of type $T$ may happen at the same time too frequently, which makes defining such martingale and infering finer probabilistic statements difficult. In order to address that, we define some auxillary processes. 

Frst, we define $\widetilde{T}_n,$ which will count only a subset of trees of type $T$ created. Let $\widetilde{T}_0=0.$ The time a tree of type $T$ is generated by vertex addition in $G_n,$ if no tree of type $T$ is removed, we flip a coin with success probability $\frac{1}{\sqrt{n}}$ and let $\widetilde{T}_{n+1}=\widetilde{T}_{n}+1.$ For the negative change in $\widetilde{T}_n,$ since the maximum negative change in $T_n$ occurs if all the neighbors of the new vertex is attached to a distinct $T,$ we can set its probability $m\frac{\widetilde{T}_n|T|}{n}.$ We clearly have $\widetilde{T}_n \leq T_n.$

Let $p_{\widetilde{T}}(G_n)$ stand for the probability that a tree of type $T$ is created and none of the type $T$ is removed from the graph. Now define a level more abstract process where we overlook the possibility of vertex addition to remove a tree of type $T.$ Let us define it as
\begin{equation*} 
\begin{aligned}
   &   
  M_{n+1} = M_n + \begin{cases}   
       1,&  \text{with prob. } \frac{p_T}{\sqrt{n}},\\[5pt]
          -m,&  \text{with prob. } \frac{M_n|T|}{2m(n+1)},        \\[5pt]
       0,&  \tn{ otherwise},
    \end{cases}
\end{aligned}
\end{equation*}
which satisfies $T_n' \leq M_n$ when coupled properly. We check the conditions of Lemma\ref{lil} for $M_n$. We take $s=m|T|.$ The variance $s_n^2$ is of order $n^{2s+\frac{1}{2}},$ and the martingale differences satisfy
\[|X_i| \leq C_1 n^s \leq \frac{n^{s+\frac{1}{4}}}{\log n} \quad \tn{a.s.}\]
Therefore, we have 
\[\limsup_n Z_n \leq C_2 n^{s+1/4} \log n,\]
which implies $ M_n \leq C\sqrt{n}$ a.s. for some $C>0$, as its mean is of order $\sqrt{n}.$ Thus, $\widetilde{T}_n\leq M_n \leq C \sqrt{n}.$

The consequence of the result above is that, since the probability of attaching $n$ to any subset of vertices of size $C \sqrt{n}$ goes to zero as $n\rightarrow \infty,$ the limiting probability for tree generation excluding the trees generated by $\widetilde{T}_n,$ will be the same. Hence, we can find $N>0$ such that $n>N$ implies $p_{\widetilde{T}}(G_n) > \frac{p_T}{2\sqrt{n}}.$ Then we define $\widetilde{M}_n$ by taking $\widetilde{M}_n=0$ for $n\leq N$ and setting
\begin{equation*} 
\begin{aligned}
   &   
  \widetilde{M}_{n+1} =\widetilde{M}_n + \begin{cases}   
       1,&  \text{with prob. } \frac{p^T}{2\sqrt{n}},\\[5pt]
          -m,&  \text{with prob. } \frac{\widetilde{M}_n|T|}{2m(n+1)},        \\[5pt]
       0,&  \tn{ otherwise},
    \end{cases}
\end{aligned}
\end{equation*}
We have $\widetilde{M}_n \leq \widetilde{T}_n \leq T_n$ when coupled properly. Now we apply the law of the iterated logarithm again to have \[\liminf_n \widetilde{Z}_n \geq C_1' n^{s+1/4} \log n,\]
which implies $T_n \geq \widetilde{M}_n \geq C'\sqrt{n}$ a.s.

Now, suppose the chosen representative isomorphism class for the logical class cannot be obtained as a tree with the root $n,$ but can be obtained from an isomorphism class of the latter type. Let us call that type $T'.$ We can define a similar martingale, this time taking generating probabilities as transition probabilities from $T$ to $T'.$ Let us call that probability $\delta_{T,T'}$ and let 
\[p_{T'}(n) := \frac{\delta_{T',T} T_n}{n}.\]

Then, since $T_n$ is almost surely of order $n^{1/2}$ as shown above, the same argument and Lemma \ref{lil} apply to $T'.$ Similarly, we can argue for the other classes that can only be obtained from the former and will end up having each isomorphism class of trees with diverging order of copies in the attachment graph.
 
\bbox

\subsection{Cycle-trees} \label{cycles}

 We will consider the cycle-trees, defined in Section \ref{evolution}, with cycles of length at most $r=2^k$ where $k$ is the fixed length of the logical sentences. Observe that what determines the logical class of a cycle is the combinations of the logical classes of the trees attached to the $l$ vertices around it. So we reduce our analysis to the trees stemming from the vertices attached to cycles of fixed length $l$ where $2 \leq l \leq 2^{k}-1.$ In case of multi-cycles, we consider the union of such trees. We fix the cycle length throughout the section. We will consider both the first occurence of those trees, that is the moment a cycle is created, and also the evolution of them.

Let $\mathbf{v}=\{v_1,\ldots,v_l\}$ be a set of vertices in $G_n$ such that $v_i\leq v_{i+1}$ and $v_l=n$. Let us define
\begin{equation} \label{hyp}
H_l^n:=\{(v_1,\ldots,v_{l-1},n) \in \{1,2,\ldots,n\}^l \, : \, v_1 < v_2 < \ldots <v_{l-1} <n\}.
\end{equation}

Next observe that any cyclic permutation of them defines a cycle that can be generated at the $n$th stage of the attachment graph, and it is unique up to taking inverses. We let 
\[\Pi_l= \tn{C}_l / \sim  \tn{ where } \sigma \sim \sigma^{-1} \]
denote the set of all cycles on $l$ elements, in other words the equivalence class of the set of cyclic permutations of length $l$ where a permutation and its inverse is identified. There are $\frac{(l-1)!}{2}$ such cycles for $l>2$.

%consider two different labellings of  We will have two different labellings of the vertices around the cycle. First,

Given $\mathbf{v} \in H_l$ and $\bar{\sigma} \in \Pi_l,$ we first order the edges around the cycle in the way that the first and the last edge have $n$ as one of their endpoints. Then for the $i$th edge let $n_i$ and $k_i$ be its endpoints such that $k_i < n_i.$ Let $\mathbf{n}=(n_1,n_2,\ldots,n_l),$  where $n_1=n_l=n$ and $\mathbf{k}=(k_1,\ldots,k_l).$ Obviously $\mathbf{k},\mathbf{n} \subseteq \mathbf{v}.$ Now we look at the probability that the prefigured cycle is generated at the $n$th stage, which is 

\begin{align*} 
\mathbf{P}(\tn{C}_l(\mathbf{v},\sigma))&= \mathbf{E}\left[ \prod_{i=1}^{l} \left( \sum_{j=1}^m \left( \alpha_{n_i}(j) \frac{1}{n_{i}-1} + (1-\alpha_{n_i}(j))\frac{D_{n_{i}}(k_i,j)}{2m(n_{i}-2)+j-1}\right)+O(n_i^{-2})\right)\right] 
\end{align*}
where the error term $\cal{O}(n_i^{-2})$ is for the correlations among $\{D_{n_i}(k_i,j) \, : \, j=1,\ldots,m\}$. We can further simplify the expression by\eqref{alpha} and combine the error terms to get 

\begin{align}\label{cycgen}
\mathbf{P}(\tn{C}_l(\mathbf{v},\sigma))&=m^l\mathbf{E}\left[ \prod_{j=1}^{l} \left( \alpha \frac{1}{n_{i}} + (1-\alpha)\frac{D_{n_{i}}(k_i)}{2mn_{i}} + O(n_i^{-2})\right)\right] \notag \\
&= \left(\frac{1-\alpha}{2}\right)^l \mathbf{E} \left[ \prod_{j=1}^{l} \frac{D_{n_{i}}(k_i)+2mu+O(n_i^{-1})}{n_{i}} \right].
\end{align}
The error term $\cal{O}(n_i^{-2})$ can be bounded below and above by universal constants independent of $\mathbf{n}$ and $\mathbf{k}.$

First, we will bound the probability above for a given cycle. Then, we will sum over all possible cycles to show that the probability of a new bounded cycle to be created is of order $n^{-1}.$ Finally, we will consider the trees of bounded depth around those cycles, which were defined as cycle-trees in Section \ref{evolution}, see \eqref{cycneigh}. Our aim is to show the existence of a limit for the generation probability of each logical class of cycle-trees. 

\subsubsection{Degree distributions}

To study the degree distribution in the expression \eqref{cycgen}, we consider its expansion in \eqref{degreedist}, which consists of products of $\beta$-random variables when conditioned on $\cal{F},$ the $\sigma$-algebra generated by $\{\psi_i\}_{i=1}^{\infty}$. We will bound the moments of the $\beta$-distribution to make inferences on the degree distributions. First note the symmetry of the $\beta$-random variables:
\[X \sim \beta(a,b) \tn{ if and only if } 1-X \sim \beta(b,a).\]

\begin{lemma} \label{prodbound}
Let $u$ and $\chi$ be defined as in \eqref{uchi} and $s\geq1.$ We have
\begin{equation*}
\left(\frac{n_1+2u-1+\frac{1}{2m}}{n_2+2u-1+\frac{1}{2m}}\right)^{s\chi} \leq \mathbf{E} \left[  \prod_{i=n_1+1}^{n_2} (1-\psi_i)^s \right] \leq \left(\frac{n_1+2u+\frac{s}{2m}}{n_2+2u+\frac{s}{2m}}\right)^{s\chi}
\end{equation*}

\end{lemma}

\pf
We first estimate the moments of 
\[\psi_i \sim \beta(m+2mu, (2i-3)m+2mu(i-1)),\]
or, equivalently of
\[1-\psi_i \sim \beta((2i-3)m+2mu(i-1), m+2mu).\]
We note that the moments of $X=\beta(a,b)$ is given by
\begin{equation}\label{betamom}
\mathbf{E}[X^s]=\prod_{i=0}^{s-1}\frac{a+i}{a+b+i} \tn{ for }s\geq 1.
\end{equation}
So, for $i\geq 2,$ we have
\begin{align*}
\mathbf{E}\left[(1-\psi_i)^s\right] &= \prod_{r=0}^{s-1} \frac{(2i-3)m+2mu(i-1)+r}{(2i-2)m+2mui+r} \\
&= \prod_{r=0}^{s-1} \left( 1-\frac{m+2mu}{(2i-2)m+2mui+r} \right) \\
& = \prod_{r=0}^{s-1} \left(1- \frac{\chi}{i+ 2u+\frac{r}{2m}-1}\right).
\end{align*}

%For the lower bound, it will suffice to consider the case where $i=\alpha n.$ 
Note the following inequalities:
\begin{equation} \label{expineq}
\left(1-\frac{1}{A} \right)^{\chi} \leq 1-\frac{\chi}{A} \leq e^{-\frac{\chi}{A}}.
\end{equation}
By \eqref{expineq} and an elementary bound on the logarithmic function, we have
\begin{align*}
\mathbf{E} \left[ \prod_{i=n_1+1}^{n_2} (1-\psi_i)^s \right] &\leq \tn{exp}\left( -s\chi \sum_{i=n_1+1}^{n_2+1} \frac{1}{i+2u-1+\frac{s}{2m}}\right) \\
&\leq \tn{exp}\left( -s\chi \log \left( \frac{n_2+2u+\frac{s}{2m}}{n_1+2u+\frac{s}{2m}}\right)\right) \\
& \leq  \left( \frac{n_1+2u+\frac{s}{2m}}{n_2+2u+\frac{s}{2m}}\right)^{s\chi}.
\end{align*} 
On the other hand, we have the lower bound
\begin{align*}
\mathbf{E} \left[  \prod_{i=n_1+1}^{n_2} (1-\psi_i)^s \right] &\geq    \prod_{i=n_1+1}^{n_2} \left(1-\frac{\chi}{i+2u-1+\frac{1}{2m}}\right)^s \\
& \geq   \prod_{i=n_1+1}^{n_2} \left(1-\frac{1}{i+2u-1+\frac{1}{2m}}\right)^{s\chi} \\
& =   \prod_{i=n_1+1}^{n_2} \left(\frac{i+2u+\frac{1}{2m}}{i+2u-1+\frac{1}{2m}}\right)^{s\chi} \\
& = \left(\frac{n_1+2u-1+\frac{1}{2m}}{n_2+2u-1+\frac{1}{2m}}\right)^{s\chi}.
\end{align*} 

\bbox

%Observe that
%\[\mathbf{E}[\psi_k^2(1-\psi_k)]=\mathbf{E}[(1-\psi_k)^2[1-2(1-\psi_k)+(1-\psi_k)^2]],\]
%so we do not need to bother about those intersections. 

\begin{lemma} \label{uppd}
For $n_1 > n_2$ and a fixed vertex $k$ in $G_n,$ we have
\[\mathbf{E}\left[D_{n_1}(k) \, D_{n_2}(k)\right] \leq \frac{\chi^2 m(m+1)}{(1-\chi)^2} \left(\frac{(n_1+A)(n_2+A)}{(k-1+\alpha)^{2}} \right)^{1-\chi}.\]
\end{lemma}
\pf 
 Recall that $\cal{F}$ is the sigma field generated by $\{\psi_i\}.$ From \eqref{chidef}, \eqref{chi} and \eqref{degreedist}, we obtain the upper bound
\begin{equation*}  
\begin{split}
\mathbf{E}[D_{n_1}(k) D_{n_2}(k) \, | \, \cal{F}] &= \sum_{a_1=k+1}^{n_1} \sum_{a_2=k+1}^{n_2} \, \sum_{i_1=1}^m \, \sum_{i_2=1}^m \mathbf{E} \left[\chi_{k_1,i_1}^{(a_1)}  \chi_{k_2,i_2}^{(a_2)} \, | \, \cal{F} \right]  \\
& \leq \sum_{a_1=k+1}^{n_1}  \sum_{i=1}^m \mathbf{E} \left[\chi_{k_1,i_1}^{(a_1)} \, | \, \cal{F} \right] \sum_{a_1=k+1}^{n_2}  \sum_{i=1}^m  \mathbf{E}\left[\chi_{k_2,i_2}^{(a_2)} \, | \, \cal{F} \right] + \sum_{a=k+1}^{n_1}  \sum_{i=1}^m \mathbf{E} \left[\chi_{k,i}^{(a)} \, | \, \cal{F} \right] \\
& \leq \mathbf{E}[D_{n_1}(k) \, | \, \cal{F}] \, ( \mathbf{E}[D_{n_2}(k) +1 \, | \, \cal{F}]).
\end{split}
\end{equation*}
This will imply that, conditioned on $\cal{F},$ their correlation is insignificant. To study the expected value of the upper bound obtained above, we will write the product of conditional in terms of $\beta$-random variables. First, let us define 
\begin{equation} \label{shifted}
D'_n(k)=D_n(k)-m
\end{equation}
to simplify the calculations. We have
\begin{equation} \label{betacorr}
\mathbf{E}[D'_{n_1}(k) \, | \, \cal{F}] \, ( \mathbf{E}[D'_{n_2}(k) \, | \, \cal{F}])  = m^2\sum_{l_1=k+1}^{n_1} \sum_{l_2=k+1}^{n_2} \psi_k^2 \prod_{i=k+1}^{l_1} (1-\psi_i)^2 \prod_{j=\min\{l_1,l_2\}+1}^{\max\{l_1,l_2\}} (1-\psi_i)
\end{equation}
Let $A=2u+\frac{1}{2m}.$ By the independence of $\{\psi_i \, : \, i=1,\ldots, n\},$ Lemma \eqref{prodbound} and \eqref{betamom}, we have
\[\mathbf{E}\left[ \mathbf{E}[D'_{n_1}(k) \, | \, \cal{F}] \, ( \mathbf{E}[D'_{n_2}(k) +1  \, | \, \cal{F}])\right] \leq \frac{4\chi^2m(m+1)}{(k-1+\alpha)^2} \sum_{l_1=k+1}^{n_1} \sum_{l_2=k+1}^{l_1}  \left(\frac{k+A}{l_1+A}\right)^{2\chi} \left(\frac{l_1+A}{l_2+A}\right)^{\chi}.  \]
Then, we apply the upper bound of the following
\begin{equation} \label{intapp}
\frac{(b)^{1-\chi}-(a-1)^{1-\chi}}{1-\chi} = \int_{a}^{b+1} \left(\frac{1}{x}\right)^{\chi} dx \leq \sum_{i=a}^b \left( \frac{1}{i}\right)^{\chi} \leq \int_{a-1}^{b} \left(\frac{1}{x}\right)^{\chi} dx= \frac{(b+1)^{1-\chi}-a^{1-\chi}}{1-\chi}
\end{equation}
twice to get the result.

\bbox 

\subsubsection{Cycle generation probabilities}

\begin{lemma} \label{cycub}
Let $\mathbf{v}=(v_1,\ldots,v_l)$ be an arbitrary set of vertices in $G_n$ except $v_l=n,$ and $\sigma$ be one of their cyclic configurations. We have the probability of generating a cycle over them
\[\mathbf{P}(\tn{C}_l(\mathbf{v},\sigma)) \leq C(\alpha,l) n^{-2\chi} v_1^{2\chi-2}\prod_{i=2}^{l-1} v_i^{-1} \]
where $C(\alpha,l)$ depends only on the graph parameter and the length of the cycle.
\end{lemma}
\pf 
We first bound the equation \eqref{cycgen} for any given $\mathbf{v}$ and $\sigma \in \Pi_l$. We can disregard the shifts by constants and the error terms in \eqref{cycgen} for the upper bound. Let $\mathbf{n}$ and $\mathbf{k}$ be the associated strings as in the begining of the section. Suppose there are $s$ pairs of vertices
\[\mathbf{l}=\{k_{1,1},\ldots,k_{s,1},k_{1,2},\ldots,k_{s,2}\} \subseteq \mathbf{k}\] 
such that $l_j:=k_{j,1}=k_{j,2},$ that is to say both $n_{j,1}$ and $n_{j,2}$ connect to $l_j.$ Let $L \subseteq [l]$ be the index set of $\mathbf{l}.$ By successively applying \eqref{difftarget}, we have 
\[\mathbf{E}\left[ \prod_{i=1}^l\frac{D'_{n_{i}}(a_{i})}{n_i} \right] \leq \prod_{j=1}^l (n_i)^{-1} \prod_{j=1}^s  \mathbf{E}\left[D_{n_1}(k) \, D_{n_2}(k)\right] \prod_{j \notin L} \mathbf{E}[D_{n_{j}}(k_j). ] \]
Then we write the above expression in terms of $\beta$-random variables and use Lemma \ref{uppd} to have
\begin{align*}
\mathbf{E}\left[ \prod_{i=1}^l\frac{D'_{n_{i}}(k_{i})}{n_i} \right]  \leq & \prod_{i=1}^l (n_i)^{-1} \prod_{j=1}^s  
\frac{4\chi^2m(m+1)}{(1-\chi)^2} \left(\frac{(n_{j,1}+A)(n_{j,2}+A)}{(l_j-1+\alpha)^{2}} \right)^{1-\chi}
 \\
& \times \prod_{j \notin L} \frac{\chi m}{1-\chi} \left( \frac{n_j+A}{k_j-1+\alpha} \right)^{1-\chi} \\
 \leq & m^{l-s}(m+1)^s 2^{-(l-2s)}\left(\frac{\chi}{1-\chi}\right)^l \prod_{i=1}^l (n_i)^{-1} \left(\frac{n_i+A}{k_i-1+\alpha} \right)^{1-\chi}\\
 \leq &  C(\alpha,l)  \prod_{j=1}^l \left( \frac{n_i}{k_i} \right)^{1-\chi} n_i^{-1}
\end{align*}
where $C(\alpha,l)$ is a constant independent of $\mathbf{n}$ and $\mathbf{k}.$ Therefore, 
\begin{align}\label{degreeUB}
\mathbf{P}(\tn{C}_l(\mathbf{v}, \sigma)) &\leq C'(\alpha,l) \prod_{i=1}^l n_i^{-1} \sum_{S \subseteq [l]} A^{|S|} \prod_{s \in S} C_s \left( \frac{n_i}{k_i} \right)^{1-\chi} \notag \\
&\leq D(\alpha,l) \prod_{i=1}^l \left( \frac{n_i}{k_i} \right)^{1-\chi} n_i^{-1}.
\end{align}

In the second part of the proof, we consider all combinations of $v_1 < \dots < v_{l-1}<v_l=n$ that give a cycle and find the one that maximizes the probability \eqref{degreeUB}. To do so, we use the idea in \cite{BHJL23} employed for uniform attachment graphs. 

For any cycle $\bar{\sigma} \in \Pi_l$, let $\sigma$ be one of the two cyclic permutations associated with it. We define a one-sided degree vector, meaning that we only consider the edges from a vertex with a higher indexed to a lower indexed one, $d_{\bar{\sigma}}=(d_{\bar{\sigma}}(1),\ldots,d_{\bar{\sigma}}(l))$ entrywise by indicator functions as
\begin{equation} \label{outgoing}
d_{\bar{\sigma}}(i):= \mathbf{1}_{\{\sigma(i) > \sigma(i+1)\}}+\mathbf{1}_{\{\sigma(i) > \sigma(i-1)\}}=\mathbf{1}_{\{\sigma^{-1}(i) > \sigma^{-1}(i+1)\}}+\mathbf{1}_{\{\sigma^{-1}(i) > \sigma^{-1}(i+1)\}}.
\end{equation}
 It is evidently well-defined. Observe that each entry is precisely the number of edges in the cycle created in the graph at the time the vertex corresponding to that entry is added. Equivalently, it is the number of outgoing edges when the attachment graph is considered as a directed graph, where the directions are simply assigned from the vertex created to the vertices it is attached to. This expression is shown to be maximized for $(d(1),\ldots,d(l))=(0,1,\ldots,1,2)$ in Lemma 6.3 of \cite{BHJL23} for uniform attachment graphs, i.e., when $\alpha=1,$
\begin{equation} \label{unifatt}
\max_{\bar{\sigma} \in \Pi_l } \prod_{i=1}^{l} v_i^{-d_{\bar{\sigma}}(i)}= v_l^{-2}\prod_{i=2}^{l-1} v_i^{-1}.
\end{equation}
The maximizer is indeed the class of the identity permutation. 

For the general case, where $\alpha > 0$, recallling that $\frac{1}{2} \leq \chi \leq 1,$ the upper bound  in \eqref{degreeUB} can be written as 
\begin{equation} \label{product}
\prod_{i=1}^l n_i^{-\chi} k_i^{\chi-1}=\prod_{i=1}^l v_i^{-\chi d_{\bar{\sigma}}(i)+ (2-d_{\bar{\sigma}}(i))(\chi-1)}.
\end{equation} 
Similarly, we have  
\begin{equation} \label{prefatt}
\max_{\bar{\sigma} \in \Pi_l} \prod_{i=1}^l v_i^{-\chi d_{\bar{\sigma}}(i)+ (2-d_{\bar{\sigma}}(i))(\chi-1)}= v_1^{2\chi-2}v_l^{-2\chi} \prod_{i=2}^{l-1} v_i^{-1}.
\end{equation}
The maximum is again achieved for the class of the identity permutation.

\bbox

\begin{proposition} \label{cycsup}
Let $\alpha>0$ and $\rho_l(n)$ be the probability that a cycle of length $l$ is created in $G_n.$ We have $\sup_{n} n \rho_l(n) < \infty.$  
\end{proposition}
\pf We will use
\[\rho_l(n) \leq \sum_{v \in H_l} \sum_{\sigma \in \Pi_l} \mathbf{P}(\tn{C}_l(\mathbf{v},\sigma)).\]
to bound the probability. For $\alpha=1,$ we have 
\[\sum_{\mathbf{v} \in V_l}\mathbf{P}(\tn{C}_l(\mathbf{v}, \tn{id})) \leq   \sum_{v_1< \ldots < v_{l-1}<n} n^{-2}\prod_{i=2}^{l-1} v_i^{-1} =  \sum_{v_2< \ldots < v_{l-1}<n} n^{-2}\prod_{i=3}^{l-1} v_i^{-1} = \cdots = n^{-1}.  \]
Recalling that $v_l=n,$ using \eqref{unifatt} and \eqref{prefatt}, we have
\begin{align*}
\sum_{\mathbf{v} \in V_l} \mathbf{P}(\tn{C}_l(\mathbf{v},\tn{id})) 
& \leq   D(\alpha,l) \sum_{v_1< \ldots < v_{l-1}<n} n^{-2\chi} v_1^{2\chi-2}\prod_{i=2}^{l-1} v_i^{-1}  \\
&=  \frac{ D(\alpha,l)}{2\chi-1} \sum_{v_2< \ldots < v_{l-1}<n} n^{-2\chi} v_2^{2\chi-2}\prod_{i=3}^{l-1} v_i^{-1} \\
&\cdots = \left( \frac{ C(\alpha,l)}{(2\chi-1)^{l-3}} \right) n^{-1}.
\end{align*}

Therefore,
\begin{align*}
\rho_l(n) & = \sum_{\sigma \in \Pi_l} \sum_{\mathbf{v} \in V_l}  \mathbf{P}(\tn{C}_l(\mathbf{v}, \sigma)) \\
& \leq  \frac{(l-1)!}{2} \sum_{\mathbf{v} \in V_l} \mathbf{P}(\tn{C}_l(\mathbf{v},\tn{id})) \\
& = D'(\alpha,l) n^{-1}
\end{align*}
for some constant $D'(\alpha,l),$ which finishes the proof.

\bbox

%\subsubsection{Integral and the use of the lower bound}

Finally, we show the main result of this section. For any ordered set of isomorphism classes of rooted trees of a fixed depth $r$, we show that the probability of creating a cycle whose $r$-neighborhood consists of such trees has a scaled limit. 

\begin{lemma} \label{cycconv}
Let $U=(T_1,\ldots,T_l)$ be a fixed isomorphism class of cycle-trees of length $l,$ and $\rho_l^U(n)$ be the probability that a cycle-tree belonging to $U$ is created in $G_n.$ We have
\[\lim_{n \rightarrow \infty} n\rho^U_l(n) \in (0,\infty).\]
\end{lemma}
\pf 
We first show the existence of the limit. For a given set of vertices $\mathbf{v}=\{v_1 < \ldots< v_l\},$ let us define the event
\[B(\mathbf{v}):=\{B_r(v_i)=T_i \, \tn{ for all }i\}\]
Then for a cycle $\sigma \in \Pi_l,$ we define 
\begin{equation}\label{event}
 U(\mathbf{v},\sigma):= \tn{C}_l( \mathbf{v},\sigma) \cap B(\mathbf{v}).
 \end{equation} 

For the cycle formation, we also need to consider the degree distribution of a vertex by the times the other vertices are added. So we define the vector random variable on the ordered pairs. For all $1 \leq i<j\leq l,$ we let
\[D_{\mathbf{v}}(i,j):=D_{v_j}(v_i).\]

Next, we define a function on 
\[\cal{H}_l:=\{(x_1,\ldots,x_{l-1},1) \in [0,1]^l \, : \, x_1 \leq x_2 \leq \ldots \leq x_{l-1} \leq 1\},\]
which is just the limiting compact space of $H_l^n$ defined in \eqref{hyp}. For 
$$x = (x_1,\ldots,x_{l-1}) \in \cal{H}_l,$$
let
\[\mathbf{v}_n(x) := \left(\ceil{n x_1},\ldots, \ceil{n x_{l-1}} \right), \]  
and
\[F_n(x):=n^l \sum_{\sigma \in \Pi_l} \mathbf{P}(U(\mathbf{v}_n(x),\sigma)).\]
Observe that
\begin{equation} \label{cycintegral}
\int_{\cal{H}_l} F_n(x) dx =  n\sum_{\mathbf{v} \in H_l^n}  \sum_{\sigma \in \Pi_l}  \mathbf{P}(U(\mathbf{v},\sigma))  = n\rho^U_l(n). 
\end{equation}

We will show that the integral converges as $n \rightarrow \infty$. We will use the dominated convergence theorem, so we start with showing the pointwise convergence of $F_n(x)$. 

Let us write 
\[\mathbf{P}(U(\mathbf{v},\sigma))=\mathbf{E}_{D_\mathbf{v}} \mathbf{E}\left[\mathbf{1}_{U(\mathbf{v},\sigma)} \, | \, D_{\mathbf{v}} \right]\]
First, we have the convergence of $D_{\mathbf{v}_n(x)}$ to a vector random variable as 
\[\big| \ceil{n x_i} < u < \ceil{n x_j} \, : \, u \tn{ is attached to } \ceil{n x_i} \big|\]
has a limiting distribution which can be obtained from \eqref{pol1}. 
Note that the weak convergence, as stated in Lemma \ref{coupling}, connects the sequential model to the P\'{o}lya graph. Conditioned on $D_{\mathbf{v}_n(x)},$ the cycle generation probability in \eqref{cycgen} is deterministic. Secondly, the probability of the second event in \eqref{event} conditioned on both $D_{\mathbf{v}_n(x)}$ and the generated cycle has a limit by Lemma \ref{coupling} since we consider trees of a fixed depth at fixed roots. Now, if we replace $\mathbf{v}$ by $\mathbf{v}_n(x)$ in \eqref{cycgen}, we have $\prod_{i=1}^l n_i= \prod_{i=1}^l \ceil{n x_{l-1}}$, which implies that the expectation in the expression multiplied by $n^l$ is finite as $n\rightarrow \infty.$ Putting them together, we conclude that $n^{l}\mathbf{P}(U(\mathbf{v}_n(x),\sigma))$ is convergent. Therefore, for every $x,$ we can define $F(x)$ as 
 \[ F(x) = \lim_n F_n(x).\]
The second step is to find an integrable function $G$ that satisfies $F_n(x) \leq G(x)$ for all $x \in \cal{H}_l.$ By Lemma \ref{cycub}, we have
\begin{align*}
 \mathbf{P}(\tn{C}_l(\mathbf{v}_n(x), \sigma))& \leq C(\alpha,l) \ceil{nx_1}^{2\chi-2}n^{-2\chi} \prod_{i=2}^{l-1} \ceil{nx_i}^{-1} \\
& \leq C(\alpha,l) n^{-l} x_1^{2\chi-2} \prod_{i=2}^{l-1} x_i^{-1}.
\end{align*} 
So, if take 
\[G(x)= \frac{(l-1)!}{2}C(\alpha,l) \, x_1^{2\chi-2}  \prod_{i=2}^{l-1} x_i^{-1} \cdot \mathbf{1}_{\{x_1 < x_2 <\ldots<x_{l-1}\}}, \]
we have $F_n(x) \leq G(x)$ since $U(\mathbf{v}_n(x),\sigma) \subseteq \tn{C}_l(\mathbf{v}_n(x),\sigma)$. If $\chi > \frac{1}{2}$ ($\alpha>0$), then  
\begin{align*}
\int_{\cal{H}_l} G(x) d x & = C' \int_{0}^{1} \int_{0}^{x_{l-1}} \cdots \int_{0}^{x_2} x_1^{2\chi-2}  \prod_{i=2}^{l-1} x_i^{-1} d\prod_{i=1}^{l-1} x_i \\
&= \frac{C'}{2\chi-1}  \int_{0}^{1} \int_{0}^{x_{l-1}} \cdots \int_{0}^{x_3} x_2^{2\chi-2}  \prod_{i=2}^{l-1} x_i^{-1} d\prod_{i=2}^{l-1} x_i \\
& \cdots = C''< \infty
\end{align*}
Therefore, by the dominated convergence theorem,
\[\lim_{n\rightarrow \infty}\int_{\cal{H}_l} F_n(x) d x \] 
exists for $\alpha >0$. Comparing it to \eqref{cycintegral}, we obtain the desired result except the case that uniform attachment has zero probability. \\

For $\alpha=0,$ the expression cannot be uniformly bounded by cycle generation probabilities for the good reason that they are of asymptotically larger order than the probabilities in the other cases. However, if an isomorphism class $U$ is fixed, they have the same order. We will study other aspects of this model separately in Section \ref{sabg}. 

For a given $\sigma \in \Pi_l,$ assume that $v_l=n$ is to be connected to $v_j = \ceil{ n x_j }$ for some $j$ according to $\sigma$ Let $d_j$ be the degree of the root of $T_j,$ and $\delta:=d_j^{-2}.$ Now observe that in the upper bound of Lemma \ref{cycub}, we computed the attachment probability using the expected value of $\mathbf{E}[D_n(v_1)].$ However, if we bound the degree of $v_j,$ 
\[D_n(v_j \, ; \, B_r(\ceil{n x_j})=T_j) \leq d_j. \] 
Secondly, observe that $\alpha=0$ is the special case that the upper bound \eqref{product} does not depend on $\sigma,$ since $-\chi=1-\chi=-1/2$ in the exponent. Therefore, having that
\[U(\mathbf{v}, \sigma) \subseteq \tn{C}_l(\mathbf{v}_n(x), \sigma) \cap \{B_r(\ceil{n x_j})=T_j)\},\]
if we replace $\sqrt{n/v_j}$ by $d_j$ for $x_j \leq \delta,$  we get
\begin{align*}
\mathbf{P}(U(\mathbf{v}, \sigma)) \leq  \mathbf{1}_{\{x_j \leq \delta\}} \cdot \ceil{n x_j}^{-1/2} n^{-3/2} \prod_{i\neq j,n} \ceil{n x_i}^{-1} +  \mathbf{1}_{\{x_j > \delta\}} \cdot \prod_{i=1}^l \ceil{n x_i}^{-1}. 
\end{align*}
Then we let
\[ G_0(x) := \begin{cases} 
     x_j^{-1/2}  \displaystyle\prod_{i\neq j}  x_i^{-1} & x_j \leq \delta \\ 
      \displaystyle\prod_{i=1}^{l-1}  x_i^{-1} & x_j > \delta, 
   \end{cases}
\]
We can show that 
\[\int_{H_l} G_0(x) dx < C(m,U),\]
which allows us to apply the dominated convergence theorem. The existence part is concluded.\\

For the positivity of the integral, first, observe that
\[\lim_{n \rightarrow \infty} U(\mathbf{v}_n(0,x_2,x_3,\ldots,x_l),\sigma)= 0\]  
since $\lim\mathbf{P}(B_r(\ceil{n x_1})=U_1)=0,$ as the degree of the first vertex goes to infinity. However, if we restrict the function to a compact subset of $H_l^n$ bounded away from zero, such as
\[K_l^n:=\{v\in H_l^n \, : \, v_1 \geq \ceil{\delta n} \}\]
for some $\delta>0.$ Let $\cal{K}_l$ be its extension to real vector space. Then, we have
\begin{equation} \label{kompakt}
Q_U(x):=\lim_{n \rightarrow \infty} \mathbf{P}\left( B(\mathbf{v}_n(x)) \right)>0
\end{equation}
by Lemma \ref{coupling}, since the probability converges to the sampling probability of the associated P\'{o}lya point-tree. Now, since $Q_U(x)$ in \eqref{kompakt} is continuous on $\cal{H}_l$ by the same lemma, it achieves a minimum in the compact set$\cal{K}_l$, let us call it $Q_U^{*}.$

Now, we have
\begin{align} \label{complim}
\lim_{n\rightarrow \infty}\int_{\cal{K}_l} F_n(x) d x & \geq Q_U^{*}\lim_{n\rightarrow \infty} n \sum_{v \in K_l^n }  \sum_{\sigma \in \Pi_l} \mathbf{P}\left(\tn{C}_l(\mathbf{v},\sigma) \, | \, B(\mathbf{v})\right) \notag \\ &   \geq Q_U^{*}\lim_{n\rightarrow \infty} n \sum_{v \in K_l^n }   \mathbf{P}\left(\tn{C}_l(\mathbf{v},\tn{id}) \, | \, B(\mathbf{v})\right) 
\end{align}
Then, since the attachment probabilities are at least $\frac{1}{n}$ regardless of the degree distribution, the case $\alpha=1$ provides a lower bound. We have 
\begin{equation} \label{loglb}
\sum_{v \in K_l^n } \mathbf{P}(\tn{C}_l(\mathbf{v},\tn{id})  \, | \, B(\mathbf{v}) ) \geq \sum_{v \in K_l^n } \mathbf{P}(\tn{C}_l(\mathbf{v},\tn{id}) \, ; \alpha=1)  =   \sum_{v \in K_l^n } n^{-2} \prod_{i=2}^{l-1} v_i^{-1}.
\end{equation}
Given that $v_l=n,$ we have
\[\sum_{v \in K_l^n } \prod_{i=2}^{l} v_i^{-1} \longrightarrow \int_{\delta}^1 \int_{\delta}^{x_{l-1}}\cdots \int_{\delta}^{x_2} \prod_{i=2}^l x_i^{-1} \ dx >0 \quad \tn{ as } \, n\rightarrow \infty.\]
This shows that the limit on the right in \eqref{complim} is positive, which implies the what is to be shown.

\bbox

\subsection{Rare occurence of multiple cycles} \label{rare}

In \cite{BHJL23}, it is shown that any leafless multicycle isomorphism class has finite distribution in the uniform attachment graph, which implies that the probability of creating two or more cycles of bounded depths is summable there. We will show the summability for the sequential model of the preferential attachment graph. 

That will allow us to condition on the high probability event that there will be no more multicycle neighborhoods of bounded depth is created for large enough number of vertices. See Section \ref{evolution} for the definition of multicycles. This fact was also essential in \cite{MZ22} for bounded degree uniform attachment graphs. The proof in \cite{BHJL23} is adaptable for the uniform attachment case, where the only difference is that we pin a vertex as in the previous section. For the preferential attachment case, we need further work.

First, we extend Proposition \ref{cycsup} to a fixed vertex other than $n.$
\begin{lemma} \label{cyclelemma}
Let $v$ be a fixed vertex in $G_n.$ We have the probability that there exists a cycle of length $l$ with one of its vertices is $v$ bounded above as 
\[\mathbf{P}(\tn{C}_l \, ; \, v \in \tn{C}_l) \leq C v^{-1}\]
for all $n,$ and $C$ is a constant depending on $\alpha$ and $l$ only.
\end{lemma}

\pf
 Let $i(v)$ denote the rank of $v$ among the vertices of the cycle, that is to say when the vertices are listed with respect to the time they were added to the graph, $w$ is the $i(w)$th one. We have the upper bound for the probability conditioned on a fixed set of vertices of the cycle and the rank of $v$:
\[\mathbf{P}(\tn{C}_l(\mathbf{v});w \in \mathbf{v}, i(v)=s) \leq C(\alpha,l,s) v_1^{2\chi-2}\prod_{i=2}^{s-1} v_i^{-1} w^{-1} \prod_{j=s}^{l-1} v_j^{-1} v_l^{-2\chi},  \]
which can be derived from  in Lemma \ref{cycub}. Then, we have
\begin{align*}
\mathbf{P}(\tn{C}_l \, ; \, v \in \tn{C}_l , i(v)=s) &\leq C(\alpha,l,s) \sum_{v_1< \cdots<v_{s-1}}\sum_{v_s< \cdots<v_{l}} v_1^{2\chi-2}\prod_{i=2}^{s-1} v_i^{-1} v^{-1} \prod_{j=s}^{l-1} v_j^{-1} v_l^{-2\chi} \\
& \leq C'(\alpha,l,s) \sum_{v_1< \cdots<v_{s-1}}  v_1^{2\chi-2}\prod_{i=2}^{s-1} v_i^{-1} v^{-2\chi} \\
& \leq C''(\alpha,l,s) v^{-1},
\end{align*}
which would also apply either $v_1=v$ or $v_l=v.$ Then summing over $s=1,\ldots,l$ we have the desired result.  

\bbox

Let us denote the existence of a path of length $s$ between $u$ and $v$ by $u \leftrightarrow_s v.$  

\begin{lemma} \label{pathlemma}
Let $v$ be a fixed vertex in $G_n.$ We have 
\[\mathbf{P}(v \leftrightarrow_l n) \leq C \frac{(\log n)^l}{n^{\chi}}.\] 
for all $n,$ and $C$ is a constant depending only on $\alpha$ and $l$ only.
\end{lemma}

\pf Let  $i(v)$ denote the rank of $v$ as defined in the proof above. We label the fixed set of $l$ vertices $\mathbf{v}$ lying on the path from $v$ to $n$ as
\[u_1<u_2 < \cdots u_{s-1} < v < w_1 < w_2 < \cdots < w_{l-s} <n\]
provided that $i(v)=s.$ 

Observe that we can represent any path of length $l$ on a fixed set of vertices in $G_n$ as its induced graph by a permutation $\pi$ on $l$ elements, such as 
\[ v_{\pi(1)} \rightarrow v_{\pi(2)} \rightarrow \cdots \rightarrow v_{\pi(l)}.\]
Set $v_{\pi_0}=v$ and $v_{\pi_{l+1}}=n.$ Similar to \eqref{outgoing}, we define a one-sided degree sequence $d$ as
 \[d_{\pi}\left(v_{\pi(i)}\right):= \mathbf{1}_{\{\pi(i) > \pi(i+1)\}}+\mathbf{1}_{\{\pi(i) > \pi(i-1)\}},\] 
 and let $e_{\pi}(w) := 2- d_{\pi}(w).$ We will use a shorthand notation for the exponent that we will use below for calculating the probabilities,
\[\delta_{\pi}(w):= d_{\pi}(i)\chi+(2-d_{\pi}(i))(1-\chi)\]
with an exception that
\[\delta_{\pi}(v):= d_{\pi}(v)\chi+(1-d_{\pi}(v))(1-\chi).\]

% Then we take $e_i(\pi) = 2- d_i(\pi)$ for $1<i<l,$ otherwise $e_i(\pi) = 1- d_i(\pi).$

%et $u_1,\ldots,u_s$ be the vertices with smaller indexed than $v$ and $w_1,\ldots,w_{l-s}$ with larger indices. We are looking at the probability that 
Then, as in Lemma \ref{cycub}, we can bound the probability conditioned on the vertex set, the rank of $v$ and a specific permutation as
\begin{equation} \label{pathbd}
\mathbf{P}(v \leftrightarrow_l n; \mathbf{v}, i(v)=s, \pi) \leq C(\alpha,l) \left(\prod_{i=1}^{s-1} u_i^{-\delta_{\pi}(u_i)}\right) v^{-\delta_{\pi}(v)}  \left(\prod_{i=1}^{l-s} w_i^{-\delta_{\pi}(w_i)} \right) n^{-\chi}
\end{equation}

We claim that
\begin{equation} \label{delta}
\Delta_{\pi}:=\sum_{i=1}^{l-s}d_{\pi}(w_i) -  \sum_{i=1}^{l-s}e_{\pi}(w_i)\geq 0.
\end{equation}
To verify it, let us take any permutation $\pi$ of the vertices in $\mathbf{v}.$ If we remove any $u_i$ from $\mathbf{v}$ and contract the path by connecting  the neighbors of $u_i$ by an edge, then the above quantity \eqref{delta} according to the new permutation on $l-1$ elements can only get smaller. The reason is that  $u_i$ and $w_j$ being adjacent  contributes only to the second term in \eqref{delta} as $w_j > u_i$  for all $i$ and $j.$ Therefore, if we remove $u_i$ one by one for all $i,$ we will have a lower bound for the expression. But since $v<w_i<n$ for all $i,$ that bound is just zero. 

Now we sum the right-side of \eqref{pathbd} over 
\[(w_1,\ldots,w_{l-s}) \, : \, v<w_1<\cdots<w_{l-s}<n.\]
Elementary calculations reveal that the total exponents 
\[\chi \sum_{i=1}^{l-s}d_{\pi}(w_i) + (1-\chi)  \sum_{i=1}^{l-s}e_{\pi}(w_i) = (l-s) - \Delta_{\pi} \left( \frac{1}{2}-\chi \right).\]
The successive additions as in the proof of Proposition \ref{cycsup} will raise the exponent by $l-s,$ the number of variables, in addition to logarithmic terms for every vertex $w_i$ with the exponent $-\delta_{\pi}(w_i)=-1.$ Therefore, we have an upper bound 
\[\sum_{v<w_1<\cdots < w_{l-s}<n} \left(\prod_{i=1}^{l-s} w_i^{-\delta_{\pi}(w_i)} \right) \leq C_P(\alpha, l, s) (\log n)^{l-s} n^{-\chi+\Delta_{\pi}\left(\frac{1}{2} - \chi \right)}.\]

Then we sum over   
\[(u_1,\ldots,u_{s}) \, : \, u_1<u_2<\cdots<u_s<v.\]
Note that 
\[\sum_{i=1}^{l}e_{\pi}(w_i) - \sum_{i=1}^{l}d_{\pi}(w_i) = \chi,\]
because including the exponent of $n,$ the differences should add up to zero. 
Symmetrically, we have the upper bound
\[\sum_{u_1<\cdots < u_{s-1}<v} \left(\prod_{i=1}^{s-1}u_i^{-\delta_{\pi}(u_i)} \right) \leq C_P'(\alpha, l, s) (\log v)^{s-1} v^{-\delta_{\pi}(v)-\Delta_{\pi}\left(\frac{1}{2} - \chi \right)}.\]

Putting all together, we have
\[\mathbf{P}(v \leftrightarrow_l n; \mathbf{v}, i(v)=s, \pi) \leq C_P''(\alpha, l, s)(\log n)^l \left( \frac{v}{n}\right)^{\Delta_{\pi}\left(\frac{1}{2}-\chi\right)} n^{-\chi}, \]
which is maximized for any $\pi$ with $\Delta_{\pi}=0.$ We have already argued for the existence of such path. Finally, we add over all permutations and all $s=1,\ldots,l$ to conclude the proof.

\bbox

%Since there are finitely many permutations of them, an upper bound will suffice.

\begin{lemma} \label{doublecycle}
For all $\varepsilon > 0,$ there exists $N(\varepsilon)$ such that for $n>N(\varepsilon)$ the probability that a new multi-cycle is created in $G_n$ is less than $\varepsilon.$ 
\end{lemma}

\pf
By the definition of a multicycle, the graph we are interested in has no leafs, i.e., vertices of degree one. Observe that for any multicycle created by the vertex addition, we connect $n$ to two possibly identical vertices, with possibly intersecting paths of bounded lengths. Then we have those two vertices lie on a cycle or two cycles of bounded length. More than two cycles or any vertices with degree one would only add additional constraints, so would lower the probability.

In the most generic case, we have two non-identical vertices, and two separate paths with two separate cycles. Multiplying the probabilities in Lemma \ref{cyclelemma} and \ref{pathlemma}, we have 
\[P \leq C(\alpha,l,s)  \frac{1}{v_0} \frac{1}{v_0'} \frac{(\log n)^s}{n^{2\chi}}.\]
Observe that we also cover the case that the two paths leading from $n$ to two vertices intersect. First summming over all $s$ and $l,$ then over all $v_0$ and $v_0',$ from which we get additional logarithmic terms, we will arrive at a summable sequence 
\begin{equation} \label{summable}
\sum_{n} \frac{(\log n)^{5r}}{n^{2\chi}}< \infty,
\end{equation}
 where $5r$ stands for the maximum number of vertices on two cycles and a path connecting them, and $\chi >\frac{1}{2}$ \eqref{uchi}. By the Borel-Cantelli Lemma, the result follows. 

The second case is that $v_0$ and $v_0'$ lie on the same cycle. Since the path attachment probability bounds given in Lemma \ref{pathlemma} are uniform, i.e., they do not depend on $v_0,$ we can separately look at the probability that 
\[\sum_{v_0,v_1} \mathbf{P}(\tn{C}; v_0, v_1 \in \tn{C}) = \sum_{v_0} \mathbf{P}(\tn{C}; v_0 \in \tn{C}) \leq C \log n,\]
and multiply it by the former to obtain a summable bound akin to \eqref{summable}. 

Finally, suppose $v_0=v_0'.$ Then the probability that there exist two different paths leading to $w_0$ from $n$ is calculated in the same way, so is the probabilities for the cycle. Note that the cycle is allowed to contain $n$ as well. 
  
  For the uniform attachment graph, we follow the proof of Theorem 3.11 in \cite{BHJL23}, where they show the probability of interest is maximized for the case that a vertex is attached to another twice and form a two-cycle. The same exact argument holds with an exception of fixing a vertex, which will eventually give us a summable sequence of order $n^{-2}.$ The result follows.
    
 \bbox

 \subsection{Lattices of cycle-trees and fixed-rooted trees}\label{lattice}
 
In this section, we will show that the limiting distribution exists for isomorphism classes of cycle-trees of bounded depth. The generation probability of each class is asymptotically of order $n^{-1}.$ In addition to cycle-trees, we will also cover the existence of an asymptotic distribution in the case of the subtrees attached to a fixed vertex. Although the frequency of the latter is larger than the harmonic order, which depends on the fixed vertex in the preferential attachment graph, we will consider only a subset of the generated trees.

\subsubsection{Cycle-trees} We define a random process to study the evolution of cycle-trees. We are interested in the combination of the logical classes of the trees whose roots are the vertices around a given cycle. Let us denote the set of feasible isomorphism classes of cycle-trees by $\cal{U}.$ We define $U \succeq U'$ if and only if $U'$ is obtained from $U$ by adding a subtree to any vertex of $U$ of depth commensurate with the distance of its root to the cycle, so that $U'$ is still of depth $2r.$ In particular, we will consider attaching the new vertex in $G_n$ to any of the non-leaves of $U$ to obtain $U'.$ 
%let us call this process \textit{vertex addition}. 
We can grade the lattice derived above by the total number of vertices of its elements. Then it would be finite at each level although it is uncountable in total. Let
\begin{equation}\label{ufo}
\cal{U}_K:=\{ U \in \cal{U} \, : \, |U|\leq K \} \tn{ and }N_K:=|\cal{U}_K|.
\end{equation}

We will define a Markov chain on $\mathbb{Z}_{\geq 0}^{N(K)}.$ Following the strategy in \cite{MZ22}, we annul the diagonal of the transition matrix. In other words, we condition on the event that $\{X_{n+1} \neq X_n\}$ and obtain a new chain $\{Y_n\}_{n \in \mathbb{N}}.$ We will describe it, then show that it has a stationary distribution. 

% $\{\rho^{U_1},\ldots, \rho^{U_{N(K)}},\tau...\}.$

See Lemma \ref{cycconv} for the definition of $\rho^U(n)$. Note that since the cycle creation probabilities depend on the current configuration of the graph, we have an inhomogeneous Markov chain. But since there are finitely many classes under consideration and as we have shown in previous section that the cycle generation probabilities are convergent, the limiting chain will be ergodic. The formal descriptions are in the proof of the following: 

\begin{lemma} \label{lattlem}
For any isomorphism class $U$ of feasible cycle-trees, the number of its copies in the preferential attachment graph has a limiting distribution.
\end{lemma}
\pf Let us denote the probability that a subtree is attached to one of the cycle-trees of type $U$ by $\tau^U(n).$ We also consider the probability that a cycle-tree from class $U'$ passes onto one in class $U,$ call that probability $\delta(U',U).$ So, we have a Markov chain
\[Y_n=\left(Y_n^{U_1},\ldots,Y_n^{U_{N(K)}}\right)\]
with the following transition probabilities. First,
\[\mathbf{P}\left(Y_{n+1}^U=Y_n^U+1, Y_{n+1}^{U'}=Y_n^{U'} \tn{ for }U\neq U'\right)= \frac{n\rho^U(n)}{Z(Y_n)}\]
where 
\[Z(Y_n)=\sum_{U \in \cal{U}_K}Y_n^U |U| +n\sum_{U \in \cal{U}_K}\rho^U(n).\]
Then,
\[\mathbf{P}\left(Y_{n+1}^U=Y_n^U-1, Y_{n+1}^{U'}=Y_n^{U'} \tn{ for }U\neq U'\right)= \frac{n\tau^U(n)Y_n^U}{Z(Y_n)}.\]
 Observe that $\lim_{n\rightarrow \infty} n\tau(n)^U > 0$ because attaching the new vertex to any vertex in the graph has probability at least $\frac{1}{n}.$ This implies 
 \[\mathbf{P}\left(Y_{n+1}^U=Y_n^U+1, Y_{n+1}^{U'}=Y_n^{U'}-1 \tn{ and } Y_{n+1}^V=Y_n^V \tn{ for }V\neq U \tn{ or }U'\right)= \frac{\delta(U',U)Y_n^{U'}}{Z(Y_n)}.\]
 
To prove its positive recurrence, let us consider the statistic $\|Y_n\|_1=\sum Y_n^U.$ From the probabilities above that statistic is dominated by the Markov chain $\{W_n\}_{n\geq 1}$ with transition probabilities given as in Section \ref{inhomMC} where we take
\[\rho(n)= \sum_{U \in \cal{U}_K} \rho^U(n) \tn{ and } \tau(n)= \min_{U \in \cal{U}_K} \tau^U(n).\]
 Since $\cal{U}_K$ is finite, both limits exist as argued earlier and the latter is also positive. So we use Theorem \ref{conv} to identify a stationary distribution for $\{W_n\}_{n\geq 1}.$ In particular, for any $N>0,$ there exists $\varepsilon$ such that
\[\limsup_{n\rightarrow \infty} \mathbf{P}(|Y_{n}| \leq N) \geq \varepsilon,\]
which is to say there exists a finite set of states which the chain returns infinitely often. Therefore, it is positive recurrent.

The aperiodicity follows from the existence of the limit for $n\rho^U(n),$ see Lemma \ref{cyclelemma}, and of $n\tau^U(n),$ which is by the definition of the attachment probabilities in \eqref{defnt}. On the other hand, the irreducibility is implied by $\lim_{n\rightarrow \infty} n\rho^U(n) >0,$ by the same lemma, and the fact that $\lim_{n\rightarrow \infty} n\tau^U(n) >0$ as argued above. Therefore, there exists a stationary distribution for the Markov chain on the finite lattice.

\bbox

\subsubsection{Pinned trees} \label{latticepin} Now we will show how to extend the same idea to the subtrees of a fixed vertex $w$ in $G_n.$ We need some finer definitions involving neighborhoods in addition to those in Section \ref{polya}. First, we define the proper $r$-neighborhood of a vertex as
 \begin{equation*} 
\overset{\circ}{B}_r(w):= B_r(w) \setminus \{w\}.
\end{equation*}
We will consider the subtrees of $w$ excluding the cycles of bounded length containing $w.$ Let
\[T_r(w):= \overset{\circ}{B}_r(w) \setminus \left \{\bigcup B_r(\tn{C}) \, : \, w \in \tn{C} \tn{ and } |\tn{C}|\leq 2r \right \}.\] 
Now, let us define the set of subtrees in a hierarchical manner as
\begin{equation} \label{subtree}
V_w(n):=\left\{ \{T_{r-1}(w') \cap \{v \in G_n : d(v,w') < d(v,w)\}\} \, : \, w \sim w'\right\}
\end{equation}
where $w \sim w'$ means there is an edge between them, or equivalently $d(w,w')=1.$ The second set of the intersection in the definition excludes the other neighbors of $w$ and the associated subtrees from the first set. As an example, if $w\in \cal{C}_0$ as defined in \eqref{nocycle}, then
 \[T_r(w) = \bigcup_{w' \sim w} T_{r-1}(w')\]
 as an induced subgraph of $G_n.$ 
 
 Next, we consider a random set of neighbors of a given vertex. Suppose we are adding a new vertex $w'$ to the graph. According to the sequential rule \ref{defnt}, for $w<w'$ and $i=1,\ldots,m,$ we have
\[p_i(w'\rightarrow w)=\alpha_{w'}(i)+\left(1-\alpha_{w'}(i)\right)\frac{\tn{deg}_{w'}(w,i)}{2m(w'-2)+i-1},\]  
which is greater than or equal to $\frac{1}{w'},$ the smallest possible  value by definition.  Let 
\[\cal{B}_w(w') : = \tn{Ber}\left(\frac{w'}{p_i(w'\rightarrow w)}  \, \Big| \, \cal{G}_n\right)\]
where $\cal{G}_n$ is simply the $\sigma$-field generated by $G_n.$ The conditioning makes $p_i(w'\rightarrow w)$ deterministic.
Finally we define 
\begin{equation} \label{harmv}
\overline{V}_w(n)=\left\{ \{T_{r-1}(w') \cap \{v \in G_n : d(v,w') < d(v,w)\}\} \, : \, w \sim w' \tn{ and }\cal{B}_w(w')=1 \right\}.
\end{equation}
 This will be used to separate out a set of subtrees attached with harmonic rate $n^{-1}$ to a fixed vertex, so that the generation rates in the process of Section \ref{lattice} will be convergent for that particular set.
 
\begin{lemma} \label{lattlempin} Let $w$ be a fixed vertex in $G_n$ and $r$ denote a fixed depth of neighnborhoods. For any isomorphism class $T$ of trees of depth $r-1,$ the number of copies of $T$ in $\overline{V}_w(n)$ has a limiting distribution. 
\end{lemma}
\pf We replace the sets and the probabilities accordingly in the proof above. By the definition in \eqref{harmv}, the probability of attaching a new subtree to $w$ is of order $n^{-1}.$  Since the probability that a subtree is attached to a given vertex in $\overline{V}_w(n)$ depends only on the degree of the vertex scaled by $n^{-1}$, the generation and the transition probabilities are well-defined, besides the former is asymptotically convergent and positive. In addition, we have a similar lattice structure by setting the maximum size of the trees to be considered, which we denoted by $K$ above. The analogous result follows.

\bbox 
 
 \section{Positive recurrence and the proof}
 
\subsection{The infinite face of the logical classes} \label{ansiktet} Let $\cal{T}$ denote the set of all isomorphism classes of trees of a fixed depth $s$. We will define two classes of trees and show that we can restrict our attention to only of them in the limit. The uniform attachment graph differs from the rest in terms of the following distinction we will make between logical classes.

Let $\alpha \in [0,1).$ We consider the cyclic trees with all the vertices around their cycles have at least $k$ representatives of each logical class of trees at the $(s-1)$st level, which are attached to vertices as subtrees. As we fix $s,$ let $\cal{T}_{\infty}$ denote the set of isomorphism classes any such vertex belongs to. Then we call the set of isomorphism classes of the associated cycle-tree $\cal{U}_{\infty},$ and the unique logical class corresponding to it $L_{\infty}$.  If we represent a cycle-tree by trees around it, we have in fact $\cal{U}_{\infty}=(\cal{T}_{\infty},\ldots,\cal{T}_{\infty})$. 

For the uniform attachment graph, the case where $\alpha=1,$ that class will not be sufficiently containing, so we recursively define a finer class. At the basis step, we consider the trees of depth one, also known as \textit{stars}. Let us denote the isomorphism class of a star with $i$ edges by $T_i^{1}.$ We define the infinite face of trees of depth one as
\[\cal{T}_{\infty}^1 = \left \{  T_{k}^{1},T_{k+1}^{1},\ldots \right \} = \bigcup_{i=0}^{\infty} T_{k+i}^{1}.\] 
Let $\cal{S}^1_{\infty}$ be the set for the logical equivalence class of $\cal{T}_{\infty}^1,$  which is the class for all stars of size larger than $k$.

Next, we consider the trees of depth two. Let $\cal{T}^2$ denote the set of all such trees. We already have an uncountably infinite set at this stage. We define the infinite face at the second level as  
\[\cal{T}_{\infty}^2 =  \left \{  T \in \cal{T}^2 \, : \, T \tn{ contains at least }k \tn{ trees in } \cal{T}_{\infty}^1 \right \}.\] 
We define  $\cal{S}_{\infty}^2$ to be the set for the logical classes of the trees in $\cal{T}_{\infty}^2.$ 
A short note is that, for the two stages that we have covered, once a tree belong to an infinite face it remains there. At the third level and onwards, that is not the case. We will need a refined probabilistic argument to address that issue.
In general, we define  $\cal{S}_{\infty}^{i}$ to be the set for the logical equivalence classes of 
\begin{equation} \label{infinite}
\cal{T}_{\infty}^i=\{T \in \cal{T}^i \, : \, T \tn{ contains at least }k \tn{ trees from each logical class in } S_{\infty}^{i-1}  \}.
\end{equation}
Then, we let $\cal{U}_{\infty}^s=(\cal{T}_{\infty}^s,\ldots,\cal{T}_{\infty}^s)$ represent the set of isomorphism classes of the cycle-trees with trees around their cycles belonging to the infinite face at the $s$th level, and let $\cal{L}_{\infty}^s$ be the set of logical classes of $\cal{U}_{\infty}^s$. 

Finally, we define the \textit{finite logical classes} $S_f^i$ as the set for the logical classes of $\cal{T}^i \setminus \cal{T}^i_{\infty}$ and the difference $\cal{L}_f:= \cal{L} \setminus \cal{L}_{\infty}$  where $\cal{L}$ is the set of all logical classes. In case we need to refer to the depth as in the uniform attachment graph, we let $\cal{L}_f^i:= \cal{L}^i \setminus \cal{L}_{\infty}^i$ where $\cal{L}^i$ is the set of all logical classes of cycle-trees of depth $i.$

\subsection{Positive recurrence}
The main result of this section is that the number of cycle-trees that do not consist of trees belonging to the infinite face has a stationary distribution. In order to show that we keep track of the cycles generated at each stage and show them they pass onto the infinite face and tend to stay there having rates commensurate with of the cycle generation. The cases where $\alpha$  is 0 or 1 have subtleties to be addressed. 
 % so that we can show  We want to show that the number of cyclic trees in finite classes are bounded.
 
 \subsubsection{Random times} \label{randomtimes}
Let us consider the events that either a new cycle of bounded length is created or a fixed subtree is attached to an already existing cycle-tree. We assume that $n$ is large enough in order to ensure that no new multi-cycle is created, see Lemma \ref{doublecycle}. That will guarantee that those two events cannot happen at the same time.

In the following sections, we will compare the transitions of the Markov chains of interest to a random event that a new cycle is created, both of which have harmonic frequency, $n^{-1}$, to occur, except $\alpha=1$. So one can imagine flipping a biased coin to decide either to move or a create a new cycle.

Since we are looking for an upper bound for the number of cycle-trees not belonging to an infinite class, we can overestimate the probability of creating a cycle as long as we derive the desired bound. On the other hand, for already existing cycle-trees, we can underestimate the probability of attaching a vertex to them. The reason is that, since we consider those cycle-trees which will remain in the finite logical classes, the attachment will only speed up their transitions to the infinite face. We formalize this below. 

By Lemma \ref{cycconv}, the probability of creating a cycle has a limit when multiplied by $n.$ So, there exists an $N$ such that $n>N$ implies $n \rho(n) \geq r$ for some $r>0.$ While the lower bound for adding a vertex to an already existing cycle-tree is $n^{-1},$ so that we can take the conditional probabilities of the two events $\frac{r}{r+1}$ and $\frac{1}{r+1}$ respectively. 
 
The time scale for the process associated with the vertex addition will be random in the sense that it will be in comparison to the cycle generation process. For instance, among the cycles of length bounde by some $l \geq 2,$  we take the $i$th cycle created and track the evolution of its cycle-tree. We look at the statistic that the number of vertices attached to that cycle-tree by the time the $N$th cycle of bounded length is created. The following is a lower bound for that statistic: 
\[K_i=\sum_{j=1}^{N-i} \tn{Geom}\left( \frac{1}{r+1}\right)\]
where the sum is over i.i.d. geometric distributions with parameter $\frac{1}{r+1},$ that is to say $K_i$ follows a negative binomial distribution. Its left tail is not nicely bounded, yet we can still get an exponentially decaying bound for a linear difference by transforming the summands to Bernoulli distributions with $\{X_j\geq 1\}$ as the event of success. Namely, we take
 \[B_j=\tn{Ber}\left(\mathbf{P}\left(\tn{Geom}\left( \frac{1}{r+1}\right)\right)\geq 1\right)=\tn{Ber}\left(\frac{1}{r+1}\right).\]
 When coupled appropriately, $K_i$ is larger than $\sum_{j=1}^{N-i}B_j=\tn{Binom}\left(N-i,\frac{1}{r+1}\right),$ and we have an exponentially decaying bound for the latter
\begin{equation} \label{subtreecount}
\mathbf{P}\left(K_i\leq \delta\frac{N-i}{r+1}\right)  \leq \mathbf{P}\left(\tn{Binom}\left(N-i,\frac{1}{r+1}\right)\leq \delta\frac{N-i}{r+1}\right) \leq e^{-(1-\delta)^2\frac{N-i}{r+1}}
\end{equation}
by a multiplicative Chernoff bound, see \cite{HR90}.\\

%We set our time to the time of cycle creation to argue that enough many trees are added to each cyclic tree in order them to linger in the infinite face.
 
We are also interested in the convergence rate of the abovementioned process with respect to the random time. A slight abstraction will be helpful in our description. Let $\{M_n\}_{n\geq 1}$ be an inhomogeneous Markov chain on a countable state space $\Omega,$ and $\cal{L}$ be a finite set of subsets of $\Omega.$ We assume that $M_n$ converges to an ergodic process with respect to the norm in Definition \ref{matrixnorm}. Let us denote the distribution of $M_n$ over $\cal{L}$ as $\pi(M_n)$ and its stationary distribution as $\pi=(\pi_1,\ldots,\pi_{|\cal{L}|}).$ In particular, we will take the state space consisting of trees or cycle-trees as in Section \ref{lattice} and the secondary subset to be the logical classes. Let   
  \[\gamma = \min_{L \in \cal{L}} \left\{\pi_L \, : \, \pi_L>0 \right \}.\]
 We will call $L \in \cal{L}$ \textit{feasible logical class} if $\pi_L \neq 0.$ Since the process is convergent, there exists $N_{\gamma}>0$ such that $N>N_{\gamma}$ implies 
 \begin{equation} \label{convtime}
 d_{TV}(\pi(X_{N},\pi) ) \leq \frac{\gamma}{2}.
 \end{equation}
 
Now, we want to express the convergence time with respect to the progression of another process with comparable probabilities. Suppose that the two processes have relative inhomogeneous probabilities $\{r_i\}_{i\geq 1},$ for example we had $r(n)=r$ in the example above. First, we note the symmetry that, for any $N$ and $\Gamma,$ 
\begin{equation} \label{absfraction}
 \mathbf{P} \left(\sum_{i=1}^{N}\tn{Geom}\left( \frac{1}{1+r_i} \right) < \Gamma \right)  = \mathbf{P} \left(\sum_{i=1}^{\Gamma}\tn{Geom}\left( \frac{r_i}{1+r_i} \right) < N \right).
 \end{equation} 
Then, observe that we can find a deterministic $\Gamma_0$ such that 
\begin{equation} \label{fraction}
\mathbf{P} \left(\sum_{i=1}^{N_{\gamma}}\tn{Geom}\left( \frac{1}{1+r_i} \right) < \Gamma_0 \right)  \geq 1- \frac{\gamma}{4}.
\end{equation}
So, if $N(\Gamma_0)$ denotes the number of steps the principal process progressed by the time the auxillary process moved by $\Gamma_0$ steps, we have 
\begin{equation} \label{gammas}
\mathbf{P}(X_{N(\Gamma_0)} \in L) \geq \frac{\gamma}{4} \tn{ for all }L \in \cal{L}.
\end{equation}

\subsection{Stationary distribution for cycle-trees}

We will show the following result in three parts according to the parameter $\alpha$ of the model: $ \alpha \in (0,1),$ $ \alpha=1$ and $ \alpha=0.$ The order is not arbitrary, but is required by the order the proofs refer to each order.

\begin{lemma}\label{cyclelimit}
Let $k \geq 1$ be a fixed quantifier rank, defined in Section \ref{model}, and $r=\frac{3^k+1}{2}.$ For any logical equivalence class $L$ of cycle-trees, either the number of cycle-trees belonging to $L$ has a finite stationary distribution or it is almost surely greater than or equal to $k$ in $G_n$ as $n \rightarrow \infty.$
\end{lemma}

\pf The three cases are covered below; the case where $ \alpha \in (0,1)$ (the mixture of the preferential and uniform rule) in Proposition \ref{pa}, $\alpha=1$ (the uniform attachment rule) in Proposition \ref{ua} and $\alpha=0$ (the preferential rule only) in Proposition \ref{sab}. The result follows.

\bbox 

\subsubsection{The preferential attachment graph} 

\begin{proposition} \label{pa}
Let $\alpha\in(0,1)$ and $k \geq 1.$ If a logical equivalence class $L\neq L_{\infty}$, we have a stationary for the number of cycle-trees belonging to $L$ in $G_n$ as $n\rightarrow \infty.$ Otherwise, it is greater than or equal to $k$ almost surely.
\end{proposition}
%Let 
%\[S_N=|\tn{C}_l(\mathbf{v},\sigma) \notin \cal{T}_{\infty}|.\]
\pf Take $n_0>0$ large enough that after the time the $n_0$th vertex is attached, no new multi-cycle is created. There exists such $n_0$ by Lemma \ref{doublecycle}. Let $N$ denote the total number of cycle-trees generated by the time the $n$th vertex is attached, and let us denote the number of those that are not in $L_{\infty}$ by $S_{N}.$ 

Let $U_i$ be the $i$th cycle-tree created between the $n_0$th and the $n$th stage. So we have
\[\mathbf{E}[S_N]=\sum_{i=1}^N \mathbf{P}\left( U_i \notin L_{\infty} \right).\]
We want to show that 
\[\sup_{N} \mathbf{E}[S_N] < \infty.\]

Let us pick one of the $l$ vertices around the cycle of $U_i$ and call it  $w_i.$ Let $\delta>0.$  We look at the subtrees of bounded depth whose roots are attached to this vertex, see \eqref{subtree}. Now, since $\alpha <1,$ the probability of attaching a subtree to $w_i$ will increase by its degree. By the time $\ceil{\frac{\delta m(N-i)}{1-\alpha}}$ subtrees are attached, the probability of attaching a new subtree will be as large as $\frac{\delta(N-i)}{n}.$  At that time, we randomly color the subtrees with $\floor{\delta(N-i)}$ different colors, and view the evolution of the subtrees in each colored set separately, which has rate at least $n^{-1}$ for each of them. From then on, a new subtree attached to $w_i$ is randomly assigned to one of the $\floor{\delta(N-i)}$ colors.

For each color, we can construct an independent process on the lattice of trees as in Section \ref{latticepin} with the harmonic rate of attachment, defined in \eqref{harmv}. Coupled properly, each of them will give a lower bound for the number of subtrees generated in the original process. We can choose the size of the trees allowed in the process large enough, which we denoted by $K$ in Section \ref{lattice}, so that each feasible logical class has a positive probability of occurence in the stationary distribution. That is to say, there exists $K_0>0$ such that 
 \[p_L(K_0):= \mathbf{P}\left(\sum_{T \in L, \, |T|\leq K_0 }\left|M^T \right| \geq 1 \right) >0 \tn{ for all }L \in \cal{L}. \]
Let us take the minimal such probability, which exists since $\cal{L}$ is finite, 
\begin{equation}\label{minclasspf}
 \gamma_0 := \min_{L \in \cal{L}^{r-1}} p_L(K_0)  >0
\end{equation}
The stationary distribution for a fixed $K_0$ exists by Lemma \ref{lattlempin}.  There also exists $\Gamma_0>0$ such that it takes $\Gamma_0$ cycle-trees to be created for the colored sets to mix with an error of at most $\frac{\gamma_0}{2}$ in the sense of \eqref{gammas}. 
  
 Then, it follows from \eqref{subtreecount}, if we choose $\delta$ small enough, there exists a constant $0< \Delta <1$ depending on $\alpha,r,m$ and $\delta$ such that  we can bound the probability that the number of subtrees attached to $w_i$ is less than or equal to $\ceil{\frac{\delta m(N-i)}{1-\alpha}}$ as below:
\[p_{i,1} := \mathbf{P}\left(\tn{Binom}\left(N-i-\Gamma_0,\frac{1}{r+1}\right)\leq \ceil{\frac{\delta m(N-i)}{1-\alpha}}\right) \leq e^{-\Delta(N-i-\Gamma_0)}.\]

Now, given that enough many subtrees are attached, which has probability at least  $1-p_{i,1},$  by \eqref{gammas} and \eqref{minclasspf}, the probability that each class has at least $k$ representatives is bounded above by 
\begin{equation} \label{p2}
p_{i,2} :=   \sum_{l=1}^{\cal{L}}\mathbf{P}\left( \tn{Binom}\left(\floor{\delta (N-i)},\frac{\gamma_0}{2}\right)\leq k \right) \leq  |\cal{L}| \, e^{\left(-\frac{\delta (N-i) \gamma_0}{4}\right)},
\end{equation}
which follows from Hoeffding's inequality provided that
 $N - i \geq \frac{4k}{\gamma_0}.$ Denote the bounded tree with root $w_i$ by $T_i.$ The probability that $T_i$ is not in the infinite face by the time the $N$th cycle-tree is created is bounded as
\[\mathbf{P}(T_{i}(N) \notin \cal{T}_{\infty})\leq p_{i,1} + (1-p_{i,1}) p_{i,2} \leq  e^{-\frac{\Delta}{2}(N-i)}+ |\cal{L}| e^{-\frac{\delta \gamma }{4}(N-i)}\]
for
\[i \leq N - \max\left\{\frac{4k}{\gamma_0}, 2 \, \Gamma_0 \right \}.\] 

Now we consider all vertices around the cycle of $U_i$. Using the union bound, we have
\[\mathbf{P}(U_i(N) \notin L_{\infty})= l \cdot \mathbf{P}(T_i(N) \notin \cal{T_{\infty}})\leq C_1e^{-C_2(N-i)}\]
for constants $C_1$ and $C_2$ independent of $N,K$ and $i.$ At the end, we sum over all cycle-trees,
\begin{equation} 
\begin{aligned}
\mathbf{E}[S_{N}]&=\sum_{i=1}^{N}  \mathbf{P}(U_{i}(N) \notin L_{\infty}) \\
&\leq \max\left\{\frac{4k}{\gamma_0}, 2 \Gamma_0 \right \}+\sum_{i=1}^{N-\max\left\{\frac{4k}{\gamma_0}, 2 \Gamma_0 \right \}} C_1e^{-C_2(N-i)} \\
& \leq  \max\left\{\frac{4k}{\gamma_0}, 2 \Gamma_0 \right \}+  C_1' \sum_{j=1}^{N} e^{-C_2j} < C
\end{aligned}
\end{equation}
for some constant $C$ independent of $N.$

Now, we will show that there exists a stationary distribution over finite logical classes, see the end of Section \ref{ansiktet} for definitions involving logical classes. 

We define a random variable $X_{N}^K$ on $\cal{L}_f.$ Let  $X_{N}^K$ be the total number of cycle-trees of size less than or equal to $K$ and belonging to $L \in \cal{L}_f$ by the time the $N$th cycle-tree is created. Note that the following argument is independent of the choice of $K$ in \eqref{minclasspf}. We have a stationary distribution for $X_{N}^K$ as $N\rightarrow \infty$ by Lemma \ref{lattlem}, let us call it $X^{K}.$ We take $S_N^{K}=\|X_N^{K}\|_1$ and $S^{K}=\|X^{K}\|_1.$  Observe that we have the monotonicity 
\[X^{K} \leq X^{K+1}\] 
for all $K$ as $\cal{U}_K \subseteq \cal{U}_{K+1},$ which is defined in \eqref{ufo}. This itself implies that there exists a limiting distribution $X$ as $K\rightarrow \infty,$ but which is possibly infinite. To rule out that possibility, first, note that we have shown above: 
\[\mathbf{E}\left[{S^{K}}\right]\leq  \sup_N \mathbf{E} \left[S_N^K \right] \leq \sup_N\mathbf{E}[S_{N}] < \infty \]
Then, by Beppo Levi's Lemma \cite{Sc96}, there exists a finite stationary distribution over the isomorphism classes of all cycle-trees, so there is also one on each fnite logical equivalence class. 

Finally, since the number of cycles diverge and those belonging to finite classes have a stationary distribution, the infinite class has obviously more than $k$ representatives almost surely, so the second part of the proposition follows too.

\bbox
  
\subsubsection{The uniform attachment graph}

In this case where $\alpha=1$, we need a more elaborate argument. The reason is that the number of subtrees attached to a fixed vertex are of the same order, logarithmic instead of polynomial, with the number of cycle-trees created.

\begin{proposition} \label{ua}
Let $\alpha=1,$ $k\geq 1$ and $r=\frac{3^k+1}{2}.$ If a logical equivalence class $L\notin \cal{L}_{\infty}^r$, we have a stationary for the number of cycle-trees belonging to $L$ in $G_n$ as $n\rightarrow \infty.$ Otherwise, it is greater than or equal to $k$ almost surely.
\end{proposition}

%for the subtrees of a fixed rooted tree.

\pf Our reference time throughout the proof will be the number of cycle-trees generated. Let $S^{(i)}_N$ be the number of cycle-trees of depth $i$ which do not belong to $\cal{L}_{\infty}^i$ by the time the $N$th cycle-tree is created. We fix $n_0$ as in the proof of Proposition \ref{pa}. To show the existence of the limiting distribution for the cycle-trees of depth $r$, we will use a recursive argument. For all $s\leq r,$ we will consider an arbitrary vertex and the subtrees of depth $s$ attached to it, as in Section \ref{latticepin}. Then it will be inductively shown that there exists a stationary distribution over $\cal{S}_{\infty}^{s}$ for those subtrees as more and more subtrees are attached to that vertex. We will call that distribution $\lambda_s.$ We will also use an auxillary stationary distribution, which is also for the pinned subtrees, but over $\cal{S}_{f}^{s},$ which will be denoted by $\mu_s.$\\
  
\noindent \textit{Basis step}: Let $s=1.$ So, the trees we are considering are stars. For a fixed vertex, the probability of attaching the new vertex to it is asymptotically $\frac{m}{n},$ while adding a new cycle is asymptotically $\frac{\rho_1}{n}$ for some $\rho_1 > 0.$ There is a single logical class in $\cal{S}^1_{\infty}$, namely the equivalence class of stars with more than $k$ edges, hence $\lambda_1$ is a trivial distribution. Each tree will be in the unique logical class by the random time given by a sum of geometric random variables with parameter $\frac{m}{m+\rho_1}.$ It can be studied as outlined in Section \ref{randomtimes} to show
\[\sup_N\mathbf{E}\left[S^{(1)}_N\right]<\infty.\]\\

%Let $v$ be the root of a set of subtrees of depth $s.$\\

\noindent \textit{Induction:} For an arbitrary vertex $v$ and its $s$-neighborhood $T_s(v),$ the logical class of $T_s(v)$ is determined by its subtrees of depth $s-1.$ Suppose $v'$ is the root of an arbitrary chosen subtree attached to $v$. Let $V_{s-1}(K)$ denote that subtree, that is 
\[V_{s-1}(K)=T_{s-1}(v') \cup \{v'\} \setminus \{w: d(v,w) < d(v',w)\}\]
by the time $K$ new cycle-trees are created after $v'$ is connected to $v$. 

Now we state the induction hypothesis: There exists a distribution $\lambda_{s-1}$ over $\cal{L}_{\infty}^{s-1}$ such that for 
\[\gamma_{s-1}:=\min_{L \in \cal{L}_{\infty}^{s-1} } \left\{\lambda_{s-1}(L) \, : \, \pi_{s-1}(L)>0 \right \}, \]
firstly, there exists $\Gamma_{s-1}$ with
\[\mathbf{P}\left(V_{s-1}(\Gamma_{s-1}) \in L \tn{ for any } L \in \cal{L}_{\infty}^{s-1} \right) \geq \frac{\gamma_{s-1}}{2}, \]
secondly,  
\[\mathbf{P}\left(V_{s-1}(M) \notin \cal{T}_{\infty}^{s-1} \right) \leq C_{s,1}e^{-C_{s,2}M} \]
for large enough $M,$ where the constants $C_{s,1},C_{s,2}$ are independent of $M$. We want to show that those hold true for the subtrees of depth $s$ of an arbitrary vertex. 

Assume the hypothesis and take an arbitrary vertex $v.$ Let $v'$ be a neighbor of it and define the subtree $V_s(K)$ with its root $v'$ as above. Take some $\delta \in (0,1)$ and $N\geq 1.$ Consider the event that $\delta(N-i)$ subtrees of depth $s-1$ are attached to $v'$ after $v'$ itself is attached to $v,$ and $\Gamma_{s-1}$ new cycle-trees are created after the last subtree is attached. Its probability is bounded above by
\[q_{i,1}:=\mathbf{P}\left(\tn{Binom}\left(N-i-\Gamma_{s-1},\frac{1}{r+1}\right)\leq \delta (N-i)\right) \leq  e^{-\Delta(N-i-\Gamma_{s-1})}\]
for some constant $0<\Delta<1$ by \eqref{subtreecount}. We then look at the probability that each logical class has at least $k$ representatives among subtrees. By a probabilistic pigeonhole argument and the first hypothesis, it is bounded as
\begin{align*}
q_{i,2} & := \mathbf{P}\left( \tn{Multinom}(\floor{\delta (N-i)},\lambda_{s-1})(L)\leq k \tn{ for some }L\in\cal{L}^{s-1}_{\infty}\right) \\
&  \leq \sum_{l=1}^{\left|\cal{L}_{\infty}^{s-1}\right|}\mathbf{P}\left( \tn{Binom}\left(\floor{\delta (N-i)},\frac{\gamma_{s-1}}{2}\right)\leq k \right) \\
& \leq  |\cal{L}^{s-1}_{\infty}| \, \tn{exp}\left(-\frac{\delta (N-i) \gamma_{s-1}}{4}\right),
\end{align*}
again, follows from Hoeffding's inequality if $i \leq N-\frac{4k}{\gamma_{s-1}}.$ Therefore, by the definition of $\cal{T}_{\infty}^s,$ see \eqref{infinite}, we have
\begin{equation} \label{expdecay}
\mathbf{P}(V_{s}(N-i) \notin \cal{T}_{\infty}^s)\leq q_{i,1} + (1-q_{i,1}) q_{i,2} \leq C_1e^{-C_2(N-i)}
\end{equation}
provided that $i \leq N - \max\left\{\frac{4k}{\gamma_{s-1}}, 2 \, \Gamma_{s-1} \right \}$. This verifies the second hypothesis at the $s$th level.

 Now, we will argue that the first hypothesis also holds at one level higher. By the second hypothesis, mimicking the proof of Proposition \ref{pa}, we first can show that there exists finitely many subtrees of $v'$ that do not belong to an infinite class. Then, following the same proof, we can derive a stationary distribution over the subtrees pinned to vertex $v'$ and belonging to $\cal{T}_f^{s-1}$. Let us call it $\mu_{s-1}.$ Observe that $\mu_{s-1}$ is the only determinant for the logical class of  a tree belonging to $\cal{T}_{\infty}^s.$ The reason is that if a tree belongs to it $\cal{T}_{\infty}^s,$ it has more than $k$ subtrees of any infinite class at the $(s-1)$st level by definition, so only the subtrees belonging to finite classes are logically distinguishing. Therefore $\lambda_s$ over $\cal{L}_{\infty}^s$ exists. Then the implications required for the first hypothesis can be derived from the random processes we discussed in Section \ref{lattice} and Section \ref{randomtimes}. \\

\noindent \textit{The end of the proof}: Now we are ready to look at the cycle-trees. Take $s=r,$ which is the final level. We can schematize the inductive steps and the final step below as
\[\mu_0=\nu_0 \rightarrow \mu_1 \rightarrow \nu_1 \rightarrow \mu_2 \rightarrow \cdots \rightarrow \mu_r \rightarrow \pi_r \]
where $\pi_r$ is the distribution over the logical classes of cycle-trees, which is to be shown to exist.

Let $U_i$ be the $i$th cycle-tree created and $v_1,\ldots,v_l$ be the vertices around its cycle. Let us denote the associated subtrees by $V_r^{(1)},\ldots, V_r^{(l)}$ respectively, which are eventually in $\cal{T}^r_{\infty}$ by exponentially decaying rates. So, by \eqref{expdecay} and the union bound, we have
\[\mathbf{P}(U_i(N) \notin \cal{L}_{\infty}^r)= l \cdot \mathbf{P}(V_r^{(i)}(N-i) \notin \cal{T}_{\infty}^r )\leq C_1e^{-C_2(N-i)}.\]
Therefore,
\begin{equation} 
\begin{aligned}
\mathbf{E}[S_N]&=\sum_{i=1}^{N}  \mathbf{P}(U_{i}(N) \notin \cal{L}_{\infty}^r)<C
\end{aligned}
\end{equation}
as in the proof above.

Finally, we use Beppo Levi's lemma again to conclude that there exists a stationary distribution over logical equivalence classes. Then, by finiteness of the cycle-trees in finite logical classes and the existence of a stationary distribution over infinite logical classes, which are finitely many by Lemma \ref{cardlogic}, each infinite class has more than $k$ representatives almost surely, which completes the proof.  

\bbox

\subsubsection{The sequential Barab\'{a}si-Albert graph} \label{sabg}
 The final case is $\alpha=0,$ or equivalently $\chi=\frac{1}{2}.$ The difficulty here is that the number of cycles of fixed length are created with frequency of order $n^{-1}\log n,$ which makes the either argument above fail. Nevertheless, we will show that those cycles typically contain an infinite class, namely the tree stemming from the vertex with the smallest index. The frequency the remaining cycles are generated turns out to be of order $n^{-1},$ so we will be able to study them as before.

We first show a fact on the degree distribution of vertices.
\begin{lemma}\label{varbound}
Let $D_n(k)$ be the degree of the vertex $k$ at the $n$th stage of the sequential Barab\'{a}si-Albert model, i.e., the case $\alpha=0$. We have
\[\tn{Var}(D_n(k)) \leq C(m) \left(\frac{n}{k^2} +\sqrt{\frac{n}{k}} \right) \]
for some positive constant $C(m)$ depending only on $m.$
\end{lemma}

\pf We will show the result using the shifted variable in \eqref{shifted}, as the variance will remain the same. It follows from \eqref{chi}, \eqref{degreedist} and \eqref{conddegreedist} that
\[\mathbf{E}\left[D_n'^2 (k)\right] \leq \mathbf{E}\left[\mathbf{E}\left[D_n' (k) \, | \, \cal{F}\right]^2 \right] + \mathbf{E}\left[D_n' (k) \right] \] 
Then, we have 
\[ \mathbf{E}\left[D_n' (k) \, | \, \cal{F}\right]^2  = 2m^2\sum_{l_1=k+1}^{n} \sum_{l_2=k+1}^{l_1} \psi_k^2 \prod_{i=k+1}^{l_2} (1-\psi_i)^2 \prod_{j=l_2+1}^{l_1} (1-\psi_i).\]
We apply Lemma \ref{prodbound}, the integral bound \eqref{intapp} and another similar integral bound for the sums obtained as in the proof of Lemma \ref{uppd} to arrive at
\[ \mathbf{E}\left[\mathbf{E}\left[D_n' (k) \, | \, \cal{F}\right]^2 \right] \leq m^2 (k+4)\left( \frac{n}{k^2}+ \frac{4}{k^2}\right) \left( 1+\frac{1}{k-1}\right)^2 \left(1+\frac{4}{n+4}\right).   \]
 On the other hand, a lower bound for the expected value of the degree itself can be found in Section 4.3. of \cite{BBCS14}:
\[\mathbf{E}\left[D_n' (k)\right] \geq m\left( \sqrt{\frac{n}{k}} -1 \right) \left( 1- \frac{1}{k} \right).\]
The result follows from the variance formula.

\bbox

Then, using Lemma \ref{varbound}, Chebyshev's inequality gives the following.
\begin{corollary}\label{cheby}
Let $D_n(k)$ be the degree of the vertex $k$ in $G_n.$ We have
\begin{enumerate}[\normalfont(i)]
\item If $k \leq \sqrt[3]{n},$ then
\[\mathbf{P}\left(\Big| D_n(k) - m\sqrt{\frac{n}{k}} \Big| \geq k^{1/4} \sqrt{\frac{n}{k^2}}  \right) \leq C(m) \frac{1}{k^{1/2}},\]
\item If $k \geq \sqrt[3]{n},$ then
\[\mathbf{P}\left(\Big| D_n(k) - m\sqrt{\frac{n}{k}} \Big| \geq \left(\frac{n}{k}\right)^{3/8}  \right) \leq C'(m) \left(\frac{k}{n}\right)^{1/4}.\]
\end{enumerate}
\end{corollary}

In the proofs for the other two cases, we could compare the cycle generation frequency to that of uniform attachment, both of which were of order $n^{-1}.$ That is not the case here. Instead, we will use the number of subtrees attached to a fixed vertex $w$ as a reference time. Let us formally define the tree that pins other subtrees:
\begin{equation}\label{tsubtree}
T_{w}(n):= \{w\} \cup \bigcup_{ T \in V_{w}(n) } T
\end{equation}
where $V_{w}(n)$ is defined in \eqref{subtree}. We first show that there are enough many representatives of each logical class among subtrees of a given tree provided that its root has a large enough degree. 

\begin{lemma} \label{fifth}
For a vertex $w$ in $G_n,$ let $T_{w}(n)$ be defined as in \eqref{tsubtree} and $\cal{T}_{\infty}$ as in Section \ref{ansiktet}. We have 
\begin{equation*} 
\mathbf{P} \left (  T_{w}(n) \notin  \cal{T}_{\infty} \, | \, \tn{deg}(v_1)=d)    \right) \leq e^{-C(l,m)\sqrt[4]{d}}
\end{equation*}
for some constant $C(l,m)$ which depends only on $l$ and $m$. 
\end{lemma}

\pf We will identify the smallest degree fixing the exponent of the upper bound on the probability. Consider the process in Section \ref{lattice} for the set of subtrees attached to $w$ as in the proof of the first case. It has a stationary distribution for bounded sizes of subtrees as shown in that section, and has also a convergence rate $N_{\gamma}$ for $\gamma>0$ in the sense \eqref{convtime}. 
%Let us call it $N$ as there is no ambiguity. 

%Regardless of the reference time, such process

By the time $K$ subtrees attached to $w$, the probabilities of attaching a new subtree is $K/2n.$ If we consider $K$ separate processes at that time, each has at least $1/2n$ probability to proceed. In other words, the new vertex is attached to any of the subtrees of fixed depth with probability at least $1/2n.$ Since we want all of them to converge, a large enough choice, such as $K^3N_{\gamma}$ moves in total, will give an exponentially decaying bound for the convergence of each process, that is
\begin{equation}\label{r0}
r_0:=\mathbf{P}\left(\tn{Binom}\left(K^3N_{\gamma},\frac{1}{K}\right) \leq N_{\gamma} \right) \leq e^{-K N_{\gamma}}
\end{equation}
by Hoeffding's inequality. So we can take $r_i=(K+i)/K$ in the expression on the left in \eqref{absfraction} and estimate how many vertices are needed to be attached to $w$ so that all processes mix. Using Chernoff bound on the sum of independent non-identical geometric distribution, the probability that more than $\Gamma$ new vertices are needed before $K^3N_{\gamma}$ new vertices are attached to any of the $K$ subtrees is bounded above as 
\begin{equation} \label{r1}
\begin{split} 
r_1:=\mathbf{P} \left(\sum_{i=1}^{K^3N_{\gamma}}\tn{Geom}\left( \frac{K}{K+i} \right) > \Gamma \right) &\leq  \inf_{t>0} e^{-t(\Gamma-K^3N_{\gamma})} \prod_{i=1}^{K^3N_{\gamma}}\frac{\frac{K}{K+i}}{\left(1-e^t\frac{i}{K+i}\right)} \\
& \leq e^{-C(m,\gamma)K}
\end{split}
\end{equation}
if we let $e^t=1+\frac{1}{2K^3N_{\gamma}}$ and $\Gamma= 2K^4 N_{\gamma}.$
Note that the domain of the moment generating function of the geometric distribution gives another restriction on $t,$ which is $e^t < (K+i)/i$ for all $i$ in the range of the product. After finding a suitable $t,$ we used the inequalities $\log(1+x)\leq x$ for $x>-1$ and $(1+1/x)^x \leq e$ to find $\Gamma$ and derive the upper bound at the end. 
This implies that if we have the degree of $w$ larger than $2K^4 N_{\gamma},$ in other words for
\begin{equation} \label{kgamma}
K \leq \sqrt[4]{\frac{\tn{deg}(v_1)}{2N_{\gamma}}},
\end{equation}
we have
 \begin{equation} \label{r2}
  r_2:=\sum_{l=1}^{\cal{L}}\mathbf{P}\left( \tn{Binom}\left(K,\frac{\gamma}{2}\right)\leq k \right) \leq  |\cal{L}| \, e^{-\frac{ \gamma}{4}K},
\end{equation}
obtained just as in \eqref{p2}. Provided that \eqref{kgamma}, we have
\begin{align*}
\mathbf{P} \left (  T_{w}(n) \notin  \cal{T}^{\infty}   \right) \leq r_1 + (1-r_1)r_2 + (1-r_1)(1-r_2)r_3 \leq e^{-C(l,m,\gamma)\sqrt[4]{d}}.
\end{align*}
We can fix any $\gamma > 0$ at the beginning and find a large enough constant $C(l,m)$ to obtain the result.

\bbox

Now we show a distinctive fact about the sequential Barab\'{a}si-Albert model. 

\begin{lemma} \label{cyclowbd}
Let $U=(\cal{T},T_2,\ldots,T_l)$ be a cycle-tree of length $l$ where $T_i$ is a unique isomorphism class of rooted trees of a fixed depth for all $i\geq 2,$ while $\cal{T}$ stands for an arbitrary tree. Let $\rho_l^U(n)$ be the probability that a cycle-tree belonging to $U$ is created in $G_n.$ If $\alpha=0,$ then
\[\lim_{n \rightarrow \infty} \frac{n}{\log n}\rho^U_l(n) \in (0,\infty).\]
\end{lemma}

\pf The proof has the same structure with of the proof of Lemma \ref{cycconv}. For a given set of vertices $\mathbf{v}=\{v_1 < \ldots< v_l\},$ we let 
\[B^*(\mathbf{v}):=\{B_r(v_i)=T_i \, \tn{ for }i\geq 2\}.\]
Then, for a cycle $\sigma \in \Pi_l,$ we define 
\begin{equation}\label{event2}
 U(\mathbf{v},\sigma):= \tn{C}_l( \mathbf{v},\sigma) \cap B^*(\mathbf{v}).
 \end{equation} 

We define a sequence of regions converging to $\cal{H}_l$ monotonically to capture the logarithmic term in the generation probability of cycles. Let 
\[\cal{H}_l^N:=\left\{(x_1,\ldots,x_{l-1},1) \in [0,1]^l \, : \, \frac{1}{N} \leq x_1 \leq x_2 \leq \ldots \leq x_{l-1} \leq 1 \right \} \nearrow \cal{H}_l.\]
For any $x = (x_1,\ldots,x_{l-1}) \in \cal{H}_l,$ we define $\mathbf{v}_{n}(x)$ and $F_{n}(x)$ as in the proof of Lemma \ref{cycconv}. In addition, let
\begin{equation} \label{cycintegral2}
I_{N,n}:=\int_{\cal{H}_l^N} F_n(x) dx =  n\sum_{v_1 \geq \ceil{N^{-1}n}}  \sum_{\sigma \in \Pi_l}  \mathbf{P}(U(\mathbf{v},\sigma)). 
\end{equation}
Observe that 
\begin{equation} \label{rho2}
I_{n,n}=n \rho^U_l(n).
\end{equation}

Now we take the isomorphism class $U$ into account. Let $d_i$ be the degree of the root in $T_i$ in $U.$ The vertex with the expected degree $d_i$ is the solution for $v_i$ in $d_i=m\sqrt{n/v_1},$ which is $v_i=m^2n/d_i^2.$ Let $\delta_i=m(4d_i^2)^{-1}.$  For $x\leq \delta_i,$ we have
\begin{equation}\label{chebyapp}
\mathbf{P}(D_{n}(\ceil{nx}) \geq d_i  ) \leq  x^{1/4}
\end{equation}
by Corollary \ref{cheby}. Then, by Lemma \ref{cycub} and \eqref{chebyapp}, we obtain
\begin{align*}
 \mathbf{P}(U(\mathbf{v}_n(x), \tn{id}))& \leq C(l,m) \ceil{nx_1}^{-1}n^{-1} \prod_{i=2}^{l-1} \left( \ceil{nx_i}^{-1} x_i^{1/4}\mathbf{1}_{\{x_i \leq \delta_i\}} +\ceil{nx_i}^{-1}\mathbf{1}_{\{x_i >\delta_i\}} \right)\\
& \leq C'(l,m) n^{-l} x_1^{-1} \prod_{i=2}^{l-1} \left(x_i^{-3/4}\mathbf{1}_{\{x_i \leq \delta_i\}} +x_i^{-1}\mathbf{1}_{\{x_i >\delta_i\}} \right).
\end{align*} 
Taking
\[G(x)= \frac{(l-1)!}{2}C'(l,m) \, x_1^{-1} \prod_{i=2}^{l-1} \left(x_i^{-3/4}\mathbf{1}_{\{x_i \leq \delta_i\}} +x_i^{-1}\mathbf{1}_{\{x_i >\delta_i\}}. \right), \]
we have $F_n(x) \leq G(x).$ Then, we integrate $G(x)$ over the restricted region.
\begin{equation} \label{gintegral}
\begin{split}
\int_{\cal{H}_l^N} G(x) d x & = C \int_{0}^{1} \int_{0}^{x_{l-1}} \cdots \int_{1/N}^{x_2} x_1^{-1} \prod_{i=2}^{l-1} \left(x_i^{-3/4}\mathbf{1}_{\{x_i \leq \delta_i\}} +x_i^{-1}\mathbf{1}_{\{x_i >\delta_i\}} \right) d\prod_{i=1}^{l-1} x_i \\
&= C' \int_{0}^{1} \int_{0}^{x_{l-1}} \cdots \int_{0}^{x_3} (\log N + \log x_2 ) \prod_{i=2}^{l-1} \left(x_i^{-3/4}\mathbf{1}_{\{x_i \leq \delta_i\}} +x_i^{-1}\mathbf{1}_{\{x_i >\delta_i\}} \right) d\prod_{i=1}^{l-1} x_i \\
& \cdots = C_1(l,m,\delta_2,\ldots,\delta_{l-1}) \log N + C_2(l,m\delta_2,\ldots,\delta_{l-1})< \infty
\end{split}
\end{equation}
Therefore, by the dominated convergence theorem,
\begin{equation} \label{limit1}
I_N:=\lim_{n\rightarrow \infty} I_{N,n} \leq \int_{\cal{H}_l^N} G(x) d x.
\end{equation}

Finally, if we run the same computations in \eqref{gintegral} for the summation in \eqref{cycintegral2} using both parts in Corollary \ref{cheby}, we can show that
\[\frac{I_{N,n}}{\log N} \geq \frac{I_{M,n}}{\log M}, \, \tn{ which implies }\, \frac{I_{N}}{\log N} \geq \frac{I_{M}}{\log M},\]
for all $n$ and $N<M.$ Therefore,
\begin{equation} \label{limit2}
\lim_{N\rightarrow \infty} \frac{I_{N}}{\log N}
\end{equation}
exists. Putting the two limits \eqref{limit1} and \eqref{limit2} together, and expressing the diagonal sequence in terms of $\rho_l^U(n)$ by \eqref{rho2}, we can show that
\[
\lim_{n \rightarrow \infty} \frac{n}{\log n}\rho^U_l(n)
\]
exists.\\

For the positivity part, as in the proof of Lemma \ref{cycconv}, we use the uniform attachment probabilities for the edges connected to $v_2,\ldots,v_{l-1}$ to obtain a lower bound. That is
\[\sum_{\mathbf{v} \in K_l^n}\mathbf{P}(U(\mathbf{v}_n,\tn{id})) \geq \sum_{\mathbf{v} \in K_l^n} \prod_{i=2}^{l} v_{i}^{-1} n^{-1}  \sqrt{\frac{v_2}{v_1}}\sqrt{\frac{n}{v_1}}= C(l,m) \frac{\log n}{n}.\]
The result follows.

\bbox

Now we prove the positive recurrence for a particular set of cycle-trees. See \eqref{tsubtree} for the definition of $T_{v_1}(n),$ which is used below.
\begin{lemma} \label{posrec2}
If $\alpha=0$ and $C^*_n=|\{\tn{C}_l(\mathbf{v},\sigma) \, : \, v_1<\ldots < v_l \leq n \tn{ and }T_{v_1}(n) \notin \cal{T}_{\infty}\}|$ is the number of cycle-trees of fixed length $l$ in $G_n$ with the additional constraint that the neighborhood of the first vertex does not belong to $\cal{T}_{\infty},$ then $\sup_n C^*_n < \infty.$ 
\end{lemma}
\pf Any configuration of the vertices for the cycle generation has the same probability if $\alpha=0$, so we consider the one corresponding to the identity permutation and multiply it by the number of possibilities, that is
\[\sum_{\sigma \in \Pi_l} \sum_{\mathbf{v}} \mathbf{P}(\tn{C}_l(\mathbf{v},\sigma))=\frac{(l-1)!}{2} \sum_{\mathbf{v} \in V_l}\mathbf{P}(\tn{C}_l(\mathbf{v},\tn{id}).\] 

We first cut out a region with respect to the position of the first vertex. For $\beta \in (0,1),$ let 
\[\cal{C}_{l,\beta} := \bigcup_n\{ \tn{C}_l(\mathbf{v}; v_1 \geq \beta n, v_l=n)\},\]
and $\rho_{l}(n,\beta n)$ be the probability of creating a cycle in  $\tn{C}_l(\mathbf{v}; v_1 \geq \beta n, v_l=n)\}.$
Recall that
\[\mathbf{P}(\tn{C}_l(\mathbf{v})) = \frac{1}{n} \prod_{i=1}^{l-1} v_{i}^{-1}.\]
Then, we have
\begin{align*}
\rho_{l}(n,\beta n) & \leq  \sum_{\floor{\beta n} \leq v_1 <\ldots < v_{l-1} <n} \frac{1}{n} \prod_{i=1}^{l-1} v_{i}^{-1} \\
& \leq \left( n^{-1}+ o\left(n^{-2}\right) \right) \int_{\beta}^1 \int_{\beta}^{x_{l-1}} \cdots \int_{\beta}^{x_{2}} \prod_{i=1}^{l-1} v_{i}^{-1} \\
& \leq C_1(\beta,l,m)n^{-1}. 
\end{align*}
From here, we can repeat the rest of the proof in Lemma \ref{cycconv}, which will allow us to use the process in Section \ref{lattice}, using which we can mimic the case $\alpha \in (0,1)$ by taking $\rho_{l}(n,\beta n)$ as the reference time. The existence of the stationary distribution for the logical classes of cycle-trees in $\cal{C}_{l,\beta}$ will follow, and a fortiori  
\[\mathbf{E}|\cal{C}_{l,\beta}| < C_1'(\beta,l,m).\]

Therefore, we have, 
\[\mathbf{E}[C_n^*] =C_1'(\beta,l,m)+ \sum_{k=1}^n \sum_{\mathbf{v} \in H_l^k} \mathbf{1}_{ \{v_1 \leq \beta k\}} \cdot\mathbf{P}(\tn{C}_l(\mathbf{v})) \cdot \mathbf{P}(V_{v_1}(n) \notin \cal{T}_{\infty} \, | \, \tn{C}_l(\mathbf{v}) ) \]
where $H_l^k$ is defined in \eqref{hyp}, not to be confused with $\cal{H}_l^k.$ So, if we consider the sum above, by Lemma \ref{fifth} and Corollary \ref{cheby}, we obtain
\begin{align*}
&  \sum_{k=1}^n \sum_{ \mathbf{v} \in H_l^k} \mathbf{1}_{ \{v_1 \leq \beta k\}}\cdot \mathbf{P}(\tn{C}_l(\mathbf{v})) \cdot \mathbf{P}(V_{v_1}(n) \notin \cal{T}_{\infty} \, | \, \tn{C}_l(\mathbf{v}) ) \\
\leq & \sum_{k=1}^n \frac{1}{k} \sum_{\mathbf{v} \in H_l^k} \prod_{i=1}^{l-1} v_{i}^{-1} \left(e^{-C(l,m)\sqrt[10]{\frac{n}{v_1}}} + \mathbf{1}_{\{v_1 \leq \sqrt[3]{n}\}} \cdot \frac{1}{\sqrt{v_1}}+ \mathbf{1}_{\{v_1 \geq \sqrt[3]{n}   \} } \cdot \left( \frac{v_1}{n}\right)^{1/4}\right) \\
\leq & \sum_{k=1}^n \frac{1}{k} \sum_{\mathbf{v} \in H_l^k} \prod_{i=1}^{l-1} v_{i}^{-1} \left(e^{-\log{\frac{C'(l,m)n}{v_1}}} + \frac{1}{\sqrt{v_1}}+  \left( \frac{v_1}{n}\right)^{1/4}\right) \\
\leq & C_2(l,m) \sum_{k=1}^n \frac{1}{k} \sum_{\mathbf{v} \in H_l^k} \left(  n^{-1}\prod_{i=2}^{l-1} v_{i}^{-1}+ v_1^{-3/2} \prod_{i=2}^{l-1} v_{i}^{-1}+ n^{-1/4}v_1^{-3/4} \prod_{i=2}^{l-1} v_{i}^{-1} \right)  \\
\leq & C_2'(l,m).
\end{align*}

\bbox

Next, we show another distinctive, threshold fact for the case $\alpha=0.$ 

\begin{lemma}\label{v1inflimiti}
Let $U=(\cal{T}_{\infty},T_2,\ldots,T_l)$ denote the isomorphism class of cycle-trees where $T_i$ is a unique isomorphism class of rooted trees of depth $r$ for all $i\geq 2,$ while $\cal{T}_{\infty}$ represents the infinite logical class. If $\alpha=0$ and $U_n$ is the number of copies of $U$ in $G_n$ as an induced graph, then
\[\lim_{n\rightarrow \infty}\mathbf{P}(U_n\geq k)=1.\]
\end{lemma}
\pf In Lemma \ref{cyclowbd}, we showed that the generation probability of cycle-trees of type $M=(\cal{T},T_2,\ldots,T_l),$ where $\cal{T}$ was an arbitrary tree, is of order $n^{-1} \log n.$ Therefore, if we use the martingale setting in Section \ref{martaatl} by taking $M_n$ to be the number of copies of $M$ in $G_n$ as an induced graph and $p(n)=\rho_l^M(n)$, we have
\[\mathbf{E}[M_n]= \frac{\rho_l^M(n)}{s} \sim \log n\]
for some $s$ depending only on the sizes of $T_2,\ldots,T_{l-1}.$ We also have 
\[s_n^2 \sim n^{2s} \log n.\]

Checking the conditions of Lemma \ref{lil}, we can show that 
\[|X_i|\leq C n^{s} \leq  K_i \frac{s_n}{\sqrt{2 \log_2 s_n^2 }} =\frac{K n^2s \log n}{\sqrt{4s \log_e 2 \log n + \log_e 2 \log \log n}} \quad \tn{a.s.}\]
for large enough $K.$ Therefore, we have almost sure divergence of $M_n$ in logarithmic scale by Lemma \ref{lil} applied for $\liminf_{n\rightarrow \infty} Z_n,$ where $Z_n$ is the martingale associated with $M_n.$ 

Finally, we observe that $C_n^*$ in Lemma \ref{posrec2} is just the number of cycle-trees of type
\[(\cal{T} \setminus \cal{T}_{\infty},T_2,\ldots,T_l),\]
 which has a stationary distribution, therefore the almost sure divergence of $M_n$ implies the almost sure divergence of $U_n$.

\bbox

Now we are ready to conclude. Let $\cal{L}^{(1)}_{\infty}$ be the set of logical classes of the cycle trees of type $(\cal{T}_{\infty},T_2,\ldots,T_l)$ where at least one of $T_2,\ldots,T_l$ belongs to a finite logical equivalence class of trees. 

\begin{proposition} \label{sab}
Let $\alpha=1$ and $k\geq 1.$ If a logical equivalence class $L \notin \cal{L}^{(1)}_{\infty} \cup \{L_{\infty}\},$ we have a stationary for the number of cycle-trees belonging to $L$ in $G_n$ as $n\rightarrow \infty.$ Otherwise, it is greater than or equal to $k$ almost surely.
\end{proposition}
\pf The first part is by using the positive recurrence in Lemma \ref{posrec2} and Beppo Levi's lemma as before. The second part follows from Lemma \ref{v1inflimiti}.

\bbox

\subsection{Proof of Theorem \ref{mainthm}}

Lemma \ref{treelimit} and Lemma \ref{cyclelimit} combined with Lemma \ref{logic} implies that only determinants of the logical classes are cycle types, if there are enough many representatives of acyclic neighborhoods. That condition holds true by Lemma \ref{treelimit}. Then, Lemma \ref{cyclelimit} shows that any cyclic logical type either has more than $k$ representatives almost surely, or has a stationary distribution. Therefore, the cycle profile is asymptotically stable. The result follows.

\bbox

\bibliographystyle{alpha}
\bibliography{patav}

\end{document}